\documentclass[aps,prb,showpacs,12pt]{revtex4-1}
\usepackage[dvips]{graphicx}
\usepackage{bm,amsmath,amssymb}
\newcommand{\bos}{\boldsymbol}
\begin{document}

\title{Monte Carlo computation of the Vassiliev knot invariant of
  degree 2 in the integral representation}
\author{Franco Ferrari}
\email{ferrari@fermi.fiz.univ.szczecin.pl}
\author{Yani Zhao}
\email{yanizhao@fermi.fiz.univ.szczecin.pl}
\affiliation{CASA* and Institute of Physics, University of Szczecin,
  Szczecin, Poland} 
\date{\today}
\begin{abstract}
In mathematics there is a wide class of knot invariants that may be
expressed in the form of multiple line integrals computed along the
trajectory $C$ describing the spatial conformation of the knot.
In this work it is addressed the problem of evaluating  invariants of
this kind in the case in which the knot is discrete, i.~e. its
trajectory is constructed by joining together a set of segments of
constant length. Such discrete knots appear almost everywhere in
numerical simulations of systems containing one dimensional
ring-shaped objects. Examples are polymers, the vortex lines in fluids
and superfluids like helium and other quantum liquids.
Formally, the trajectory of a discrete knot is a piecewise smooth
curve characterized by sharp corners at the joints between contiguous
segments. The presence of these corners spoils the topological
invariance of the knot invariants considered here and prevents the
correct evaluation of their values. To solve this problem, a smoothing
procedure is presented, which eliminates the sharp corners and
transforms the original path $C$ into a curve that is everywhere
differentiable. The procedure is quite general and can be applied to
any discrete knot defined off or on lattice.
This smoothing algorithm is applied to the computation of the
Vassiliev knot invariant of degree 2 denoted here with the symbol
$\varrho(C)$. 
This is the simplest knot invariant that admits a definition in terms
of multiple line integrals. For a fast derivation of $\varrho(C)$, it
is used a Monte Carlo integration technique. It is shown that, after
the smoothing, the values of $\varrho(C)$ may be evaluated with an
arbitrary precision. Several algorithms for the fast computation of
the Vassiliev knot invariant of degree 2 are provided.
\end{abstract}
\maketitle
\clearpage

\section{Introduction}\label{section1}
There are many situations
in which it is necessary to distinguish the topological properties of
ring-shaped quasi one-dimensional objects. This is for instance the
case of polymers \cite{wasserman,delbruck}, vortex structures in
nematic liquid crystals \cite{vortexLC},
$^3He$ superfluid \cite{vortexHe3} and disclination lines in chiral
nematic colloids \cite{defectsNColloids}. 
In order to ascertain the type of a knot, it is possible to apply the
so-called knot invariants. These mathematical quantities, which
remain unchanged
under ambient isotopy,
are usually represented in the form of polynomials,
like for example the Alexander \cite{Alexander} or the HOMFLY
polynomials \cite{HOMFLY}. Alternatively, certain
knot invariants may be
defined in terms of multiple curvilinear integrals, in which the
integrations are 
performed along the spatial trajectory of the knot or elements of
it \cite{mci,LR,GMM}. Particularly important for  applications is the
case in which
knots are constructed by joining together at their ends a set of $N$
 segments. Discrete knots of this kind are in fact the most common
concrete realizations of knots in 
numerical simulations. Formally, a discrete knot is 
$C^0-$curve which is
piecewise smooth and is characterized by sharp corners at the joints
between contiguous segments. 
While there exist already well established
mathematical algorithms in order to compute numerically polynomial
knot invariants, see for instance the pioneering work
\cite{vologodskii}, there are not many studies concerning 
the numerical 
computation of knot invariants given in the form of multiple line
integrals for such discrete knots. 
Of course, the calculation of line integrals over discrete data is a
textbook subject \cite{reginald,davis}. Moreover,
problems in which knots are
discretized using splines, have been investigated for
example in \cite{atkinson}. However, we are faced with a somewhat
different problem, which arises due to the fact that
knot invariants 
expressed as multiple curvilinear integrals 
are  not well defined in the case of discrete knots. 
The reason of this failure is
related to the  presence of
the non-smooth corners at the joints between two contiguous
segments.
As a consequence,  the main goal pursued here is to
replace the piecewise smooth curves
representing a discrete knot with more regular ones. 
To obtain
 a smoothing algorithm that is able to perform this
replacement for general discrete knots
without destroying their topology, a
strategy has been adopted that can be briefly summarized as follows. First,
the sharp corners 
are surrounded with spheres whose radii are chosen in such a way
that they do not intersect with themselves and with other elements of
the knot. After that, 
the elements containing the corners
 inside the spheres are substituted with
arcs of smooth curves.
This  procedure 
transforms the original trajectories into $G^1-$curves\footnote{We
  recall that a $G^1-$curve is a tangent vector geometrically
continuous curve characterized by the fact that the unit tangent
vector to the curve is continuous \cite{Knott}.} without
altering their topological configurations.
As an application, the case of the Vassiliev knot invariant of
degree 2~\cite{Dunin} of a knot $C$, denoted here $\varrho(C)$, is
worked out. The
main advantages 
of choosing this invariant are its relative simplicity and the fact
that its exact values for different knots can be computed
analytically.  In this way, it is possible to perform a comparison
between  numerical  and analytical results.
After the smoothing procedure proposed in this work,  it
becomes possible 
 to calculate $\varrho(C)$ numerically
with an arbitrarily
high precision.

 Despite its
simplicity, the Vassiliev invariant of degree 2 requires the
evaluation of complicated quadruple and triple line integrals.
Having in mind concrete
applications, in which the knot invariant 
$\varrho(C)$ must be computed millions of times, see for instance
 Ref.~\cite{YZFFJSTAT2013}, we have proposed here several strategies 
to accelerate
its calculation.
It turns out that Monte Carlo integration algorithms are
 faster than traditional integration methods \cite{traditional}. For
 this reason, a
 Monte Carlo 
 integration scheme is adopted and explained in details.
Moreover, several tricks to speed up the computation
of $\varrho(C)$, that are  specific to
particular applications or situations, are presented.
Since the time for evaluating
 $\varrho(C)$ is sensitive to
the number of 
segments $N$ 
composing the knot, but not on its length, we have provided
 an algorithm to
reduce by a factor three the number of segments
without changing
the topology of the knot. 
This algorithm is valid for knots defined on a simple cubic lattice.
Secondly, 
it is shown that the 
number of points of the trajectory $C$ to be
sampled during the Monte Carlo integration procedure may be
considerably decreased 
 when the
knot invariant $\varrho(C)$ is used in order to detect topology
changes that may potentially occur after a random transformation of an
element of the trajectory of the knot.
Such random transformations, like for instance the pivot moves
\cite{pivot}, the
pull moves \cite{pull} and the BFACF moves \cite{BFACF}, are
extensively exploited in numerical simulations of 
polymer knots.

The material presented in this paper is divided as follows.
In the next Section, the Vassiliev invariant of degree 2
is defined in the case of general smooth curves. In
Section \ref{section3} 
we specialize to  general discrete knots, which are represented
as piecewise smooth curves parametrized by a continuous variable
$S\in[0,N]$. In this way, the
calculation of $\varrho(C)$ is reduced to that of multiple
 integrals over a set of variables
$S,T,U,V\in[0,N]$ and can be tackled by standard Monte Carlo techniques.
A numerical version of the so-called framing \cite{Witten}
procedure is implemented in order to  regularize 
singularities that are possibly  arising in some of the terms to be integrated.
While it is analytically proven that the sum of all these
terms is always finite, the presence of singularities in single terms
may spoil the result of the numerical integration.
A smoothing procedure is presented
in Section \ref{section4}  in order to
transform a 
general discrete knot 
 into a $G^1-$curve. This procedure allows the calculation of
$\varrho(C)$ by Monte Carlo integration techniques with an arbitrary
precision depending on the number of used samples.
In Section \ref{section5} a few methods to speed up the calculations
are discussed. Finally, the conclusions are drawn in Section
\ref{section6}. 
\section{The Vassiliev invariant of degree 2}\label{section2}
Let us consider a general knot of length $L$ in the flat three dimensional
space $\mathbb{R}^3$
spanned by a set of cartesian coordinates $\boldsymbol
x=(x^1,x^2,x^3)$. 
The space indexes  
are labeled with greek letters $\mu,\nu,\rho,\ldots=1,2,3$.
The
Alexander-Briggs notation for denoting knots is used. 
In this Section, the spatial trajectory $C$ formed in the space by the knot
is chosen to be a smooth
curve $\boldsymbol x(s):[0,L]\longrightarrow \mathbb{R}^3$
parametrized using its arc-length $0\le s\le L$.
Different points on the curve corresponding to different values of the
arc-length $s,t,u$ and  $v$ will be denoted with the symbols 
 $x^\mu(s), y^\nu(t), z^\rho(u)$ and $w^\sigma(v)$, with
$\mu,\nu,\rho,\sigma=1,2,3$. As a convention, summation over repeated
indexes is understood. Moreover,
let $\dot x^\mu(s)$  be the derivative of $x^\mu(s)$ with
respect to $s$. An analogous notation holds for $ \dot y^\nu(t),
\dot z^\rho(u)$ and $\dot w^\sigma(v)$. 
Finally,
$\epsilon_{\mu\nu\rho}$  is the completely
antisymmetric tensor uniquely defined by the condition
$\epsilon_{123}=1$. 

With the above settings, the Vassiliev knot 
invariant of degree 2 $\varrho(C)$ of a knot $C$ can be written as follows
\cite{LR,GMM,cotta}:
\begin{equation}
\varrho(C)=\varrho_{1}(C)+\varrho_{2}(C) \label{rhodef}
\end{equation}
where  $\varrho_1(C)$ and $\varrho_2(C)$ are
 two path ordered multiple line integrals
given by: 
\begin{equation}
\varrho_1(C)=\int_0^Lds \int_0^sdt\int_0^t du F_1(\boldsymbol
x(s),\boldsymbol y(t),\boldsymbol z(u);\dot{\boldsymbol x}(s),
\dot{\boldsymbol y}(t),\dot{\boldsymbol z}(u)
) \label{rhoonemi}
\end{equation}
and
\begin{equation}
\varrho_2(C)=\int_0^Lds \int_0^sdt\int_0^t du\int_0^u dv F_2
(\boldsymbol
x(s),\boldsymbol y(t),\boldsymbol z(u),\boldsymbol w(v);\dot{\boldsymbol x}(s),
\dot{\boldsymbol y}(t),\dot{\boldsymbol z}(u), \dot{\boldsymbol w}(v)
) 
\label{rhotwomi}
\end{equation}
The quantities $F_1$ and $F_2$ are defined below:
\begin{eqnarray}
-32\pi^3F_1(\boldsymbol x,\boldsymbol y,\boldsymbol z;
\dot{\boldsymbol x},\dot{\boldsymbol y},\dot {\boldsymbol z}
)&=&C_1C_2C_3\left[\dot{\boldsymbol{y}}\cdot\dot{\boldsymbol{z}}(\dot{\boldsymbol{x}}\cdot\boldsymbol{c})
+
\dot{\boldsymbol{x}}
\cdot\dot{\boldsymbol{z}}(\dot{\boldsymbol{y}}\cdot\boldsymbol{b})
+
\dot{\boldsymbol{x}}\cdot\dot{\boldsymbol{y}}(\dot{\boldsymbol{z}}\cdot\boldsymbol{a})\right]\nonumber\\
&-
&C_1C_2^2C_3\left[\dot{\boldsymbol{y}}
\cdot(\boldsymbol{a}\times\boldsymbol{b})\left(\boldsymbol{a}+
\boldsymbol{b}\frac{a}{b}\right)
\cdot(\dot{\boldsymbol{z}}\times\dot{\boldsymbol{x}})\right.\nonumber\\
&+&\left.
\dot{\boldsymbol{z}}\cdot(\boldsymbol{a}\times\boldsymbol{b})\left
(\boldsymbol{b}+\boldsymbol{a}\frac{b}{a}\right)\cdot(\dot{\boldsymbol{y}}\times\dot{\boldsymbol{x}})\right]\nonumber\\
&+&
C_1C_2\left[\dot{\boldsymbol{y}}\cdot(\boldsymbol{a}\times\boldsymbol{b})
\left(\boldsymbol{b}\frac{c-a}{b^2}+\boldsymbol{c}
\frac{a+b}{c^2}\right)\cdot(\dot{\boldsymbol{z}}\times\dot{\boldsymbol{x}})
\right.\nonumber\\
&+&\left.\dot{\boldsymbol{z}}\cdot(\boldsymbol{a}\times
\boldsymbol{b})\left(\boldsymbol{a}\frac{c-b}{a^2}-\boldsymbol{c}
\frac{a+b}{c^2}\right)\cdot(\dot{\boldsymbol{y}}
\times\dot{\boldsymbol{x}})\right]\label{F1def}
\end{eqnarray}
\begin{eqnarray}
F_2(\boldsymbol x,\boldsymbol y,\boldsymbol z,\boldsymbol w;
\dot{\boldsymbol x},\dot{\boldsymbol y},\dot {\boldsymbol z},
\dot{\boldsymbol w}
)&=&\frac{1}{8\pi^2}\left(\dot{\boldsymbol{x}\cdot}
\left(\dot{\boldsymbol{z}}\times\frac{\boldsymbol{b}}{b^3}\right)\right)
\left(\dot{\boldsymbol{y}\cdot}\left(\dot{\boldsymbol{w}}
\times\frac{\boldsymbol{c}}{c^3}\right)\right)\label{F2def}
\end{eqnarray}
In Eqs.~(\ref{F1def}) and (\ref{F2def}) we have put:
\begin{equation}
\boldsymbol a=\boldsymbol y-\boldsymbol x\qquad
\boldsymbol b=\boldsymbol z-\boldsymbol x\qquad
\boldsymbol c=\boldsymbol y-\boldsymbol z
\end{equation}
and
\begin{eqnarray}
C_1&=&\frac{2\pi}{abc}\\
C_2&=&\frac{1}{ab+a_{\mu}b_{\mu}}\label{C2}\\
C_3&=&a+b-c\label{C3}
\end{eqnarray}
Let us note that in Eq.~(\ref{C2}) we have used a convention for which
repeated indexes are summed.
It is known that the above defined $\varrho(C)$ is related
to the second 
coefficient $a_2(C)$ of the Conway polynomial of a knot $C$ through
the following 
relation \cite{GMM}: 
\begin{equation}
a_2(C)=\dfrac{1}{2}\left[\varrho(C)+\dfrac{1}{12}\right]\label{rela2rho}
\end{equation}
The coefficients of the
Conway polynomials can be
computed analytically for every knot topology.
$\varrho(C)$ is the simplest  knot invariant 
 expressed in the form of contour integrals. It is also called the
Casson knot invariant, see Ref.~\cite{polyakviro}
\section{The Vassiliev knot invariant of degree 2 for discrete
  knots}\label{section3}
In principle, analytical computations of $\varrho(C)$ are possible if the curve
$\bos x(s)$ describing the knot $C$
is given in parametric form. However,
a close expression of $\bos x(s)$  
for a knot of arbitrary shape
usually does not exist and one should pass to a discrete
representation of it. To switch from continuous functions to discrete
ones is one of the most standard problems of numerical integration.
In the present case, the situation is opposite.
For instance, in
numerical simulations involving ring-shaped objects with nontrivial
topological configurations the knots
are already discrete by construction. 
The real difficulty is rather that knot invariants expressed in the
form of multiple line integrals like $\varrho(C)$, cease to be
topological invariants if knots are discrete.
In order to restore this invariance, a
procedure 
that is able to smooth up a
discrete knot transforming it from a 
$C^0-$curve into a $G^1-$curve without destroying its topological
configuration is needed.
Such a procedure will be presented in the following.

\begin{figure}
\centering
\includegraphics[width=0.5\textwidth]{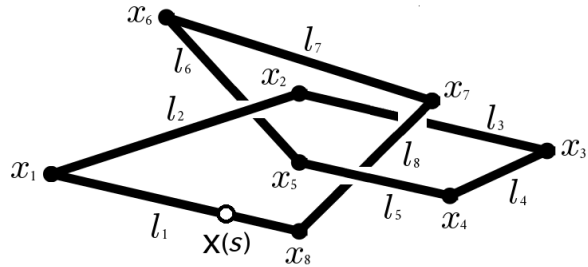}
\caption{Example of an off lattice discrete knot (a trefoil) with only eight
  sides. A generic point $\boldsymbol X(S)$ on  the trajectory
is shown.} \label{disccurv}
\end{figure}

First of all, to be concrete, let us define the discrete knot as
a set of $N$ points:
\begin{equation}
\boldsymbol x_i=\boldsymbol x(s_i)\qquad\qquad
\left\{
\begin{array}{c}
i=1,\ldots,N\\
 0< s_1<s_2\cdots <s_N=L
\end{array}
\right.\label{equ11}
\end{equation}
%
joined together by $N$ segments 
\begin{eqnarray}
\boldsymbol l_i&=&\boldsymbol x_i-\boldsymbol x_{i-1}\qquad\qquad
i=2,\ldots,N \\
\boldsymbol l_1&=&\boldsymbol x_1-\bos x_N
\end{eqnarray}
The discrete knot may be regarded as a piecewise smooth curve
$\bos X(S):[0,N]\longrightarrow \mathbb{R}$,
where
\begin{equation}
0\le S\le N\label{newpar}
\end{equation}
Explicitly, a general point located on the $i-$th segment of $\bos X(S)$
is identified by the relations:
\begin{equation}
\bos X(S)=\bos x_{i-1}+(S-[S])\bos l_i\qquad
\left\{\begin{array}{c}
i-1< S\le i\\
i=2,\ldots,N
\end{array}
\right.\label{Xsdef}
\end{equation}
and
\begin{equation}
\bos X(S)=\bos x_N+(S-[S])\bos l_1\qquad 0<S\le 1\label{Xsdefadd}
\end{equation}
In the above equations $[S]$ denotes the integer part of $S$.
The example of a curve $\bos X(S)$ with eight segments is given in
Fig.~\ref{disccurv}.
In the limit in which $N$ approaches infinity and the
lengths of the $N$ segments become vanishingly small, a continuous
representation of the knot is obtained.
If
$l_i=|\bos l_i|$ denotes the length of the $i-$th segment and
$\Lambda_N=\sum_{i=1}^Nl_i$ is the total length of the discretized
curve, then the length $L$ of the continuous knot is given by:
\begin{equation}
\lim_{\substack {N\to\infty\\ l_i\to 0,\, i=1,\ldots,N}}\Lambda_N=L
\end{equation}
At this point 
it is possible to compute the contributions
$\varrho_1(C)$
and $\varrho_2(C)$ to the Vassiliev invariant of degree 2
for a general discrete knot $C$ with trajectory $\bos X(S)$.
With the above definitions, the
same prescriptions of Eqs.~(\ref{rhoonemi})
and (\ref{rhotwomi}), which are valid for a smooth curve $\bos x(s)$,
 can be formally applied. It
is sufficient to substitute the 
smooth trajectories $\boldsymbol x(s),\boldsymbol y(t),
\boldsymbol z(u)$ and $\boldsymbol w(v)$ with their discrete analogs
 $\boldsymbol X(S),\boldsymbol Y(T),
\boldsymbol Z(U)$ and $\boldsymbol W(V)$. 
In the following, the symbols $F_1(S,T,U)$ and $F_2(S,T,U,V)$ 
will denote the integrands of
Eqs.~(\ref{rhoonemi}) 
and (\ref{rhotwomi}) in the case of a discrete knot in which
the variables $s,t,u,v$ are replaced by $S,T,U,V$.
Of course, in these equations the upper integration
boundary $L$ should be replaced by $N$.
The
derivatives $\dot{\bos X}(S),\dot{\bos Y}(T),\dot{\bos Z}(U)$ and
$\dot{\bos W}(V)$ require some more care. 
On the
$i-$th segment, away from the joints, the curve is trivially smooth
and
the computation of $\dot{\bos X}(S),\dot{\bos Y}(T),\dot{\bos Z}(U)
,\dot{\bos W}(V)$
is straightforward:
\begin{eqnarray}
\dot{\bos X}(S)&=&\bos l_i\qquad
\left\{\begin{array}{c}
i-1< S< i\\
i=2,\ldots,N
\end{array}
\right.\label{Xsderdef}\\
\dot{\bos X}(S)&=&\bos l_1\qquad 0<S< 1\label{Xsderdefadd}
\end{eqnarray}
At the
points $\bos x_1,\ldots,\bos x_N$ in which the segments join
together, instead, the curve $\bos X(S)$ ceases to be differentiable.
Still, it is possible to define 
formally the
derivatives at these points by assuming that the tangent to the
discrete
trajectory
in $\bos x_{i-1}$ is proportional to the segment $\bos l_i$.
Using this convention we obtain:
\begin{eqnarray}
\dot{\bos X}(i-1)&=&\bos l_{i}\qquad\qquad i=2,\ldots,N\\
\dot{\bos X}(N)&=&\bos l_1\label{Xsnderdef}
\end{eqnarray}
The above definition is clearly not unique. Analogously, 
we could have  
chosen
$\dot{\bos X}(i-1)=\bos l_{i-1}$, $i=2,\ldots,N$
and $\dot{\bos X}(N)=\bos l_N$.

With the prescriptions (\ref{newpar}--\ref{Xsdefadd}) 
and (\ref{Xsderdef}--\ref{Xsnderdef}) given above in order to
parametrize the discrete knot, 
the evaluation
of the two multiple line integrals appearing
in Eqs.~(\ref{rhoonemi})
and (\ref{rhotwomi}) may be performed using numerical integration
techniques like the rectangle rule method, trapezoidal rule method,
Simpson's rule method, 
Newton-Cotes method, Romberg method, Gauss method etc.~\cite{QuHe}
Let us note that
the variables $S, T$ and $U$ appearing in
$\varrho_1(C)$, see Eq.~(\ref{rhoonemi}), span a
space of 
volume 
\begin{equation}
V_1=\frac{N^3}{6}\label{vol1}
\end{equation}
 while the variables $S, T, U$ and $V$
appearing in
Eq.~(\ref{rhotwomi}) span a space of volume 
\begin{equation}
V_2=\frac{N^4}{24}\label{vol2}
\end{equation}
When $N$ is large, these volumes become too large
to be treated with
quadrature methods and
 it is more convenient
to  compute
the right hand sides of Eqs.~(\ref{rhoonemi}) and
(\ref{rhotwomi}) using a Monte Carlo approach\footnote{Nonstandard
  methods
like that based on Particle Swarm
  Optimization proposed in \cite{QuHe} 
  could probably also be applied successfully.}.
For the integral of a function 
of $m$ variables $f(\xi_1,\cdots,\xi_m)$ 
with integration boundaries like those in  Eqs.~(\ref{rhoonemi}) and
(\ref{rhotwomi}), it can be applied to this purpose
the
general formula:
\begin{eqnarray}
&&\int_{a_1}^{b_1}d\xi_1\int_{a_2}^{\xi_1}d\xi_2\cdots\int_{a_m}^{\xi_{m-1}}
d\xi_m f(\xi_1,\cdots ,\xi_m) \nonumber\\&&
\approx \dfrac{1}{n}\left[\sum_{i=1}^n f(\xi_1^{(i)},\cdots
  ,\xi_m^{(i)})(b_1-a_1)\prod_{\sigma=2}^m
  (\xi_\sigma^{(i)}-a_\sigma)\right]\label{mcf} 
\end{eqnarray}
where the $\xi_{\sigma}^{(i)}$'s, $i=1, \cdots, n$ and
$\sigma=1,\cdots,m$ denote randomly chosen variables in the range: 
\begin{equation}
\begin{array}{ccl}
[a_1,b_1]&\mbox{when}&\sigma=1\\
{[a_{\sigma},\xi_\sigma]}&\mbox{when}&\sigma=2,\ldots,m
\end{array}
\end{equation}
The naive procedure discussed above is plagued by two
systematic errors.
First of all, the discrete
knots treated here so far are not smooth at the joints between two
segments. 
On a simple cubic lattice, it is possible to verify that 
the values of $\varrho(C)$ computed for a discrete knot are always
greater than the analytical values, a fact that is certainly related
to the presence of sharp corners at these joints. This excess from the
exact value is  
roughly proportional to the number of the corners.
The second source of errors is connected with possible
singularities arising in the integrands  $F_1(S,T,U)$ and
$F_2(S,T,U,V)$ 
appearing in  Eqs.~(\ref{rhoonemi}) and
(\ref{rhotwomi}).
Of course, globally both $F_1(S,T,U)$ and $F_2(S,T,U,V)$
are regular for every value of $S,T,U$ and $V$ as it has been proved
in \cite{GMM}. 
However, the fact that $F_1(S,T,U)$ and $F_2(S,T,U,V)$ are finite
everywhere does not prevent the presence of singularities in single
terms entering in the expressions of these integrands.
Looking at Eqs.~(\ref{F1def}--\ref{C3}), it is easy to realize that
some of these terms diverge
whenever one or more of the
following conditions are met: 
\begin{eqnarray}
\boldsymbol Y(T)-\bos X(S)&=&0\label{cond1}\\
\boldsymbol Z(U)-\bos X(S)&=&0\\
\boldsymbol Y(T)-\bos Z(U)&=&0\\
|\bos Y(T)-\bos X(S)|\,|\bos Z(U)-\bos X(S)|+(\bos Y(T)-\bos
X(S))\cdot(\bos Z(U)- 
\bos X(S))&=&0\label{cond4}
\end{eqnarray}
When summed together, these singularities disappear making
$F_1(S,T,U)$ and $F_2(S,T,U,V)$ finite, but for the purposes of
numerical calculations a regularization is needed to remove them.
To 
this purpose, a suitable
 regularization  is the framing  of the trajectories described
 in \cite{Witten}.
In the present context, the framing consists in a slight deformation of the
curves $\bos X(S)$, $\bos Y(T)$, $ \bos Z(U)$ 
and $\bos W(V)$ of the kind:
\begin{eqnarray}
X^\mu(S)&\longrightarrow&X_{\epsilon_X}^\mu(S)=X^\mu(S)+\epsilon
n^\mu(S)\label{reg1}\\
Y^\nu(T)&\longrightarrow&Y_{\epsilon_Y}^\nu(T)=Y^\nu(T)+2\epsilon
n^\nu(T)\\
Z^\rho(U)&\longrightarrow&Z_{\epsilon_Z}^\rho(U)=Z^\rho(U)+3\epsilon
n^\rho(U)\\
W^\sigma(V)&\longrightarrow&W_{\epsilon_W}^\sigma(V)=W^\sigma(V)+4\epsilon
n^\sigma(V)\label{reg4} 
\end{eqnarray}
where  $n^\mu(S),n^\nu(T),n^\rho(U)$ and
$n^\sigma(V)$ denote unit vectors normal to the trajectories
$\bos X(S)$, $\bos Y(T)$, $\bos Z(U)$ and $\bos W(V)$ respectively.
$\epsilon$ 
is a very small parameter.
Clearly, the prescription provided in Eqs.~(\ref{reg1}--\ref{reg4})
is able to remove the divergences
at the locations defined in  (\ref{cond1}--\ref{cond4}).
Moreover, in the limit $\epsilon\to
0$ one recovers the exact expression of $\varrho(C)$ independently of
the choice of the normal unit vectors $n^\mu(S),n^\nu(T),n^\rho(U)$ and
$n^\sigma(V)$ as it has been proved in Ref.~\cite{GMM}.
For example, in calculations on a simple cubic lattice the framing can be
implemented by small shifts of
the trajectories $\bos X(S)$, $\bos Y(T)$, $\bos Z(U)$ and $\bos W(V)$
along the direction 
$(1,1,1)$. This is sufficient to regularize all potentially divergent
terms in  $F_1(S,T,U)$ and $F_2(S,T,U,V)$ without creating  dangerous
intersections
of the trajectories of the shifted knots that are forbidden. 
From our simulations it turns out that the results of the computations
of $\varrho(C)$ are not much sensitive 
to the values of the $\epsilon-$parameter. This is connected to the
fact that the points in which the singularity conditions of
Eqs.~(\ref{cond1}--\ref{cond4}) are satisfied represent a very small
subset of the set of all sampled points.

To eliminate the systematic error due to the presence of the sharp
corners is much more difficult. This will be the subject of
 the next Section, in which a smoothing procedure will be presented, that
transforms the curve $\bos X(S)$ 
into a curve whose first derivatives
exist and are continuous.
\section{Monte Carlo evaluation of $\Large 
\varrho(C)$ with
smoothing procedure
of dicrete knots}\label{section4}

\begin{figure}
\centering
\includegraphics[width=0.5\textwidth]{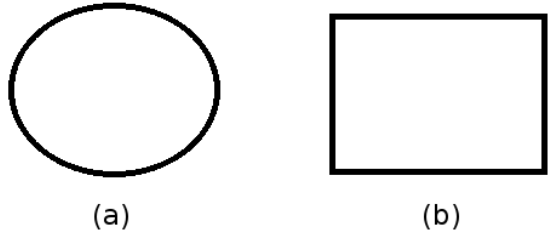}
\caption{(a) An unknot with smooth trajectory; (b) 
  An unknot defined on a simple cubic lattice.} \label{fig1}
\end{figure}

The effect of the sharp corners at the joints of the segments on the
computation of $\varrho(C)$ 
can be checked using the very simple
example
of an unknot with two
different trajectories: 
\begin{itemize}
\item A smooth circle defined by the parametric curve
$x^1(\theta)=\cos(\theta)$, $x^2(\theta)=\sin(\theta)$ and $x^3(\theta)=0$,
$\theta\in[0,2\pi]$, see
Fig.~\ref{fig1}(a).
\item A square defined on a simple cubic lattice as shown in
Fig.~\ref{fig1}(b). 
\end{itemize}
The exact value of the Vassiliev invariant of degree 2 for the unknot is
$-\frac 1{12}\sim -0.083$.
The Monte Carlo computation of $\varrho(C)$ 
gives a result that is very near to the exact one in the case of the circle:
$\varrho(C)=-0.083\pm 1.72\times 10^{-4}$.
However, for the square we obtain
$\varrho(C)=0.050\pm 1.17\times 10^{-4}$, which is far from the expected
result.
To avoid these ambiguities in the calculation of $\varrho(C)$ for
discrete knots,
a smoothing procedure for eliminating the sharp
corners will be presented.  

\begin{figure}
\centering
\includegraphics[width=0.7\textwidth]{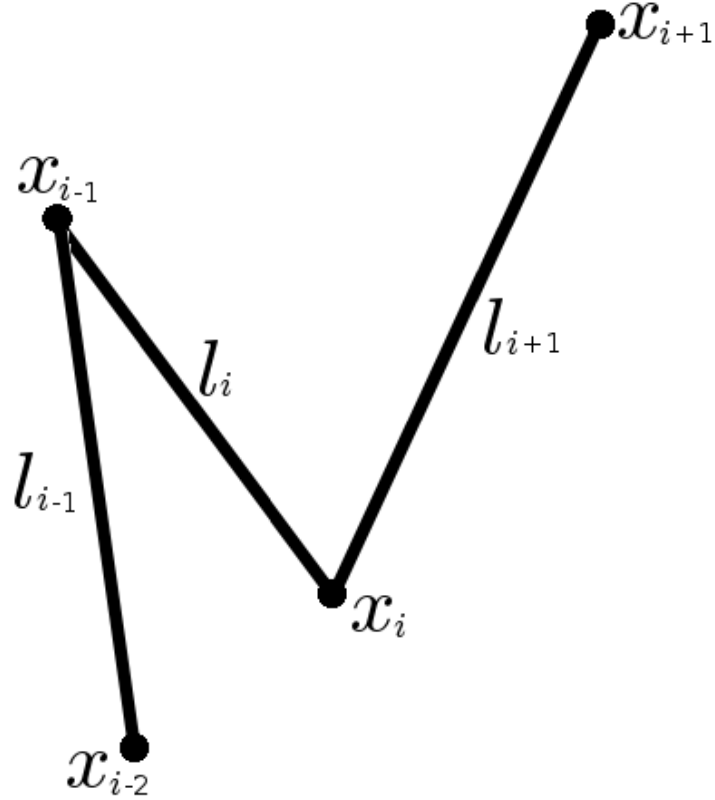}
\caption{ This figure shows the three contiguous segments $\bos l_{i-1},\bos
  l_i$ and $\bos l_{i+1}$ subtending the corners $\bos x_{i-1}$ and
  $\bos x_i$.} \label{triplet}
\end{figure}

\begin{figure}
\centering
\includegraphics[width=0.7\textwidth]{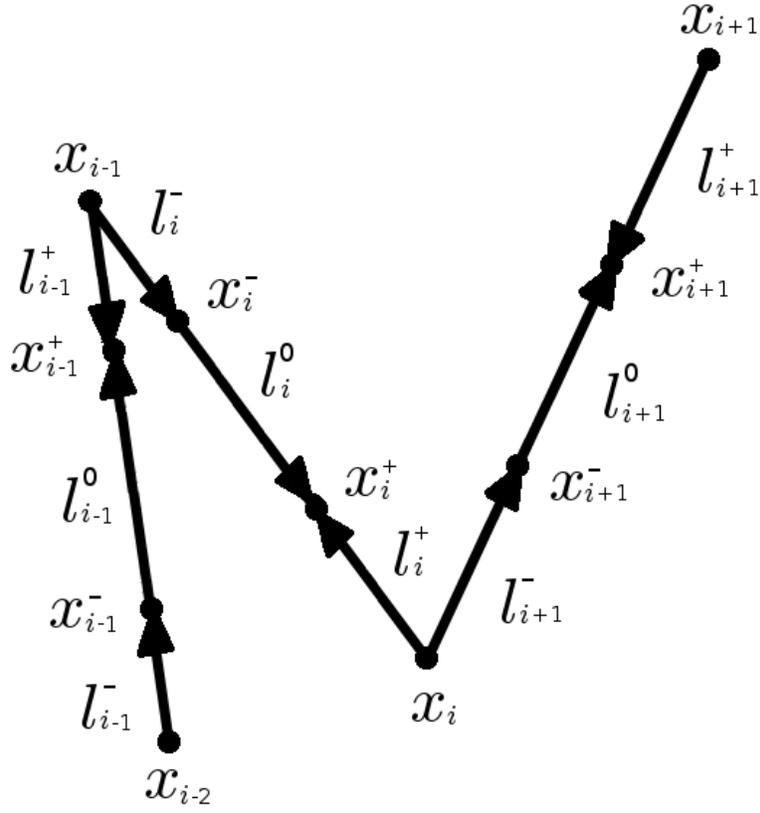}
\caption{ The segments $\bos l_{i-1},\bos l_i$ and $\bos l_{i+1}$ of
  Fig.~\ref{triplet} are split into three subsegments in such a way
  that
the corners in $\bos x_{i-1}$ and $\bos x_i$ are subtended by
segments that are not in common. In the given example, after the
splitting, the corner in $\bos x_{i-1}$ is subtended by
the subsegments $\bos l_{i-1}^+$ and $\bos l_i^-$. The corner
in $\bos x_{i}$ is subtended instead by
$\bos l_{i}^+$ and $\bos l_{i+1}^-$.} \label{splittedcorners} 
\end{figure}

The idea is to replace at each of the joints $\bos x_{i}$ 
the neighborhoods of the corners
 with  smooth arcs of curves whose ends are glued together 
in such a way that the whole knot
will be a continuous curve with its first
derivative. 
To illustrate the method, we pick up
a triplet of contiguous segments $\bos
l_{i-1},\bos l_i$ and $\bos l_{i+1}$, see Fig.~\ref{triplet}.
It is easy to realize that the segment $\bos l_i$ is shared by the two
corners 
centered at the points $\bos x_{i-1}$ and $\bos x_i$. This is not
desirable for our purposes.
To achieve the goal that each corner will be subtended by couples of segments
that are not in common with those related to other corners, we
 divide each segment $\bos l_i$, $i=1,\ldots,N$, into 
 the
three subsegments:
\begin{eqnarray}
\bos l_i^-&=&\bos x_i^--\bos x_{i-1}\\
\bos l_i^0&=&\bos x_i^+-\bos x_i^-\\
\bos l_i^+&=&\bos x_i^+-\bos x_i
\end{eqnarray}
The ends $\bos x_i^-$ 
and  $\bos x_i^+$  are 
fixed in such a way that the
lengths of $\bos l_i^-,\bos l_i^0$ and $\bos l_i^+$
are  $d_{i-1}'$, $l_i-d_{i-1}'-d_i$  and $d_i$
respectively (see
Fig.~\ref{splittedcorners}):
\begin{eqnarray}
\bos x_i^-&=&\bos x_{i-1} +\frac{\bos x_i-\bos x_{i-1}}{l_i}d_{i-1}'\\
\bos x_i^+&=&\bos x_{i} +\frac{\bos x_{i-1}-\bos x_{i}}{l_i}d_{i}
\end{eqnarray}
The values of $d_{i-1}'$ and $d_i$ will be chosen in such a way that

1) the topology of the discrete knot is not destroyed after the
smoothing procedure
and 

2) the length of none of the subsegments $\bos l_i^\pm$  and $\bos
l_i^0$ exceeds $\frac {l_i}{2}$.

\begin{figure}
\centering
\includegraphics[width=0.5\textwidth]{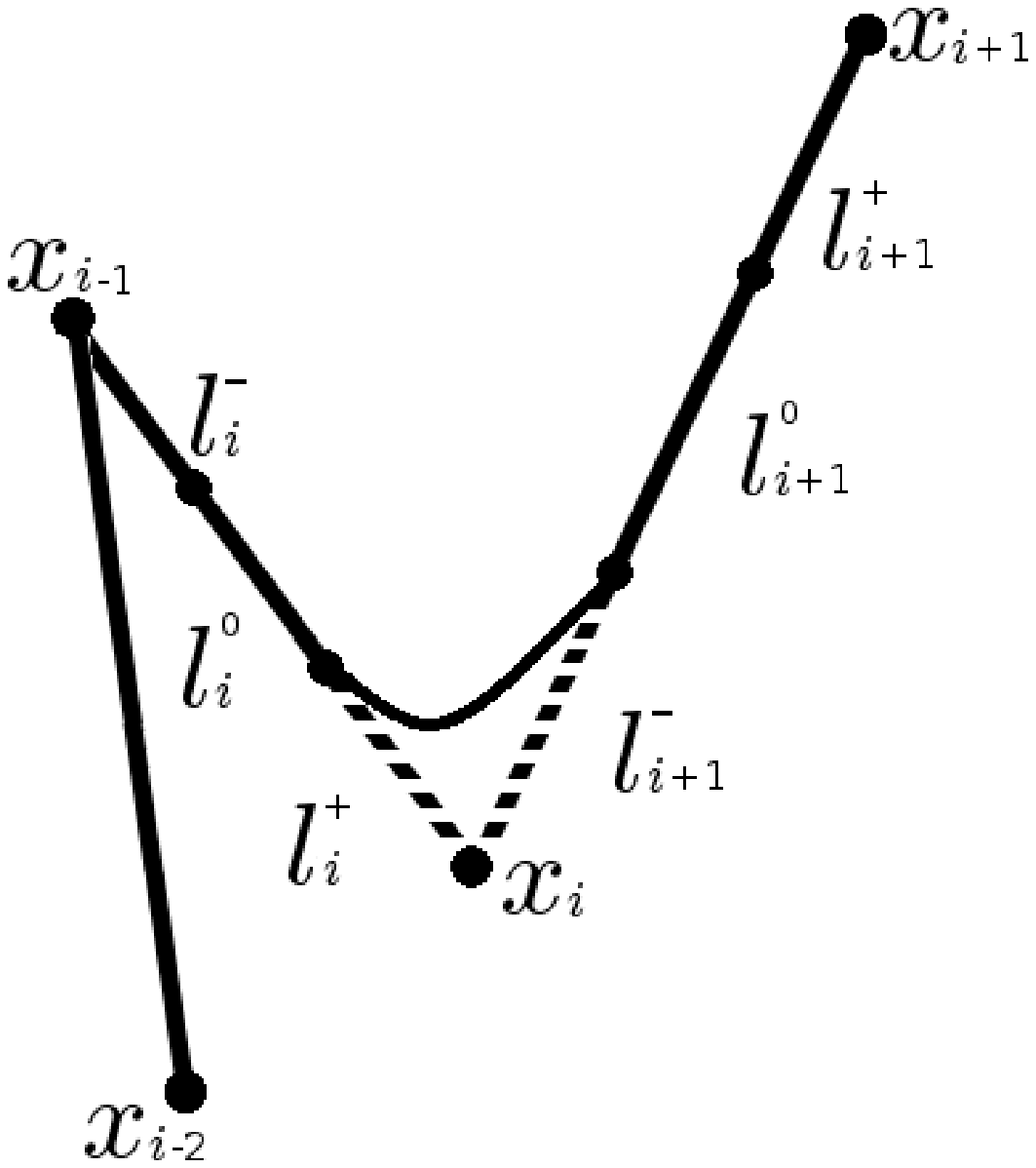}
\caption{Substitution of the sharp corner in $\bos x_i$ by an arc of
  the smooth curve defined in Eqs.~(\ref{xones}--\ref{xtwos}). It is shown that
  the subsegments $\bos l_i^+$ and $\bos l_{i+1}^-$ subtending this
  corner are replaced by a smooth trajectory. The replaced part has
  been denoted with dashed lines.} \label{offlattice} 
\end{figure}

An algorithm to determine $d_{i-1}'$ and $d_i$ will be provided later.
After performing the above splitting for $\bos l_{i-1}$, $\bos l_i$
and $\bos l_{i+1}$, the subsegments
$\bos l_{i-1}^+$ 
and $\bos l_i^-$ 
subtend the corner centered in $\bos x_{i-1}$, while the corner in
$\bos x_i$ is subtended by $\bos l_i^+$ and $\bos l_{i+1}^-$.
Thus, if all segments composing the knot are splitted in this way,
we arrive at the desired situation in which none of the segments
subtending a given corner is shared by another corner.
At this point, each corner subtended by the couples of segments
$\bos l_i^+,\bos l_{i+1}^-$ for $i=1,\ldots,N-1$ and $\bos l_N^+,\bos
l_1^-$ for $i=N$ may be substituted by arcs of smooth curves as shown
in Fig.~\ref{offlattice}.

 As an arc of a smooth curve replacing the generic corner in
$\bos x_i$ subtended by the subsegments
$\bos l_i^+$ and $\bos l_{i+1}^-$, it is possible to use
the ansatz (the superscript $+$ refers to $\bos l_i^+$, while the
superscript $-$ refers to $\bos l_{i+1}^-$):
\begin{eqnarray}
\bos
X_i^+(S)&=&-\frac{d_i}{l_i}\frac{l_{i+1}}{d_{i}'}\frac{\sin\theta_i^+(S)
}{1-\frac{1}{\sqrt{2}}+\frac{1}{\sqrt{2}}\frac{d_i}{l_i}
\frac{l_{i+1}}{d_{i}'}}\bos
  l_i^+-\frac{(\cos\theta_i^+(S)-1)}{1-\frac{1}
{\sqrt{2}}+\frac{1}{\sqrt{2}}\frac{l_i}{d_i}\frac{d_{i}'}
{l_{i+1}}}\bos
  l_{i+1}^-+\bos
l_i^++\bos x_i \label{xones}\\ 
\bos X_{i+1}^- (S)&=&-\frac{(\sin\theta_{i+1}^-(S)-1)}{1-\frac{1}{\sqrt{2}}+
\frac{1}{\sqrt{2}}\frac{d_i}{l_i}\frac{l_{i+1}}{d_{i}'}}\bos
  l_i^+-
\frac{l_i}{d_i}\frac{d_{i}'}{l_{i+1}}\frac{\cos\theta_{i+1}^-(S)}
{1-\frac{1}{\sqrt{2}}+\frac{1}{\sqrt{2}}\frac{l_i}{d_i}
\frac{d_{i}'}{l_{i+1}}}\bos
  l_{i+1}^-+\bos l_{i+1}^-+\bos x_i\label{xtwos}
\end{eqnarray}
with 
\begin{eqnarray}
\theta_i^+(S)=\left(\frac{l_i}{2d_i}(S-[S])-\frac{l_i-d_i}{2d_i}
\right)\frac{\pi}{2} &\qquad\qquad& \frac{l_i-d_i}{l_i}\le S-[S]\le 1\\
\theta_{i+1}^-(S)=\left(\frac{l_{i+1}}{2d_{i}'}(S-[S])+
\frac 12\right)\frac{\pi}{2} &\qquad\qquad& 0\le S-[S]\le 
\frac{d_{i}'}{l_{i+1}}\label{theta2}
\end{eqnarray}
Eqs.~(\ref{xones}--\ref{theta2}) are defined for
$i=1,\ldots,N-1$. Their extension 
to the corner in $\bos x_N$ is  straighforward.
It is easy to verify that after replacing the sharp corners at the
vertices 
$\bos x_i$ with
 the  arcs of curve $\bos X_i^+(S)$ and
$\bos X_{i+1}^-(S)$ a $G^1-$curve is obtained:
\begin{enumerate}
\item First of all, at the  point connecting $\bos
X_i^+(S)$
and $\bos X_{i+1}^-(S)$, occurring when $\theta_i^+(1)=\theta_{i+1}^-(0)=\frac\pi
4$, it is possible to verify that the curve obtained after the
replacement is continuous.
\item Second, both $\bos
X_i^+(S)$
and $\bos X_{i+1}^-(S)$ are differentiable and their derivatives,
which are 
continuous, do coincide.
\item Third, the unit tangent vectors computed on the subsegments
$\bos l_i^0$ coincide with the unit tangent vectors computed in the
  subsegments $\bos l_i^+$ at the point
$\bos x_i^+$
 in which these subsegments are
  connected together. The same is true in the case of the point 
$\bos x_i^-$ in  which $\bos l_i^-$ and $\bos l_i^0$ are joined.
To show that, we note that
$\bos X_i^+(S)$ maps the subsegment $\bos l_i^+$ 
into a continuous arc
of a curve with unit tangent vector $\bos t_i^+$ at the end point
$\bos x_i^+$ given by 
$\bos t_i^+=-\frac{\bos l_i^+}{l_i^+}$.  
To have a $G^1-$ curve, $\bos t_i^+$ must
coincide with the tangent $\bos t_i^0$
computed at $\bos x_i^+$, but staying   on the subsegment
$\bos l_i^0$. 
It is easy to check using the parametrization 
(\ref{Xsdef}--\ref{Xsdefadd}) of the knot on $\bos l_i^0$ that $\bos
t_i^0= \frac{\bos l_i}{l_i}$. Thus, taking into account the fact that
$\bos l_i^+$ and $\bos l_i$ are antiparallel, it is possible to
conclude that $\bos t_i^+=\bos t_i^0$ as desired.
As well, the unit tangent vector computed on the curve $\bos
X_{i+1}^-(S)$ 
at 
$\bos x_{i+1}^-$ coincides with the unit tangent vector computed on
$\bos 
l_{i+1}^0$ at the point
$\bos x_{i+1}^-$.  
\item Finally, even if this is not necessary for the present purpose, we have also checked numerically
that, for a wide 
range of the variable $x=\frac{d_il_{i+1}}{d_{i+1}l_i}$ entering the
expressions of $\bos
X_i^+(S)$
and $\bos X_{i+1}^-(S)$, more precisely
for $0.01\le x\le 100$, the distance between the point $\bos x_i$
and any of the points of the curve $\bos
X_i^+(S)$ ($\bos X_{i+1}^-(S)$) never grows beyond a fraction of $d_i$ ($d_i'$).
\end{enumerate}

\begin{figure}
\centering
\includegraphics[width=0.7\textwidth]{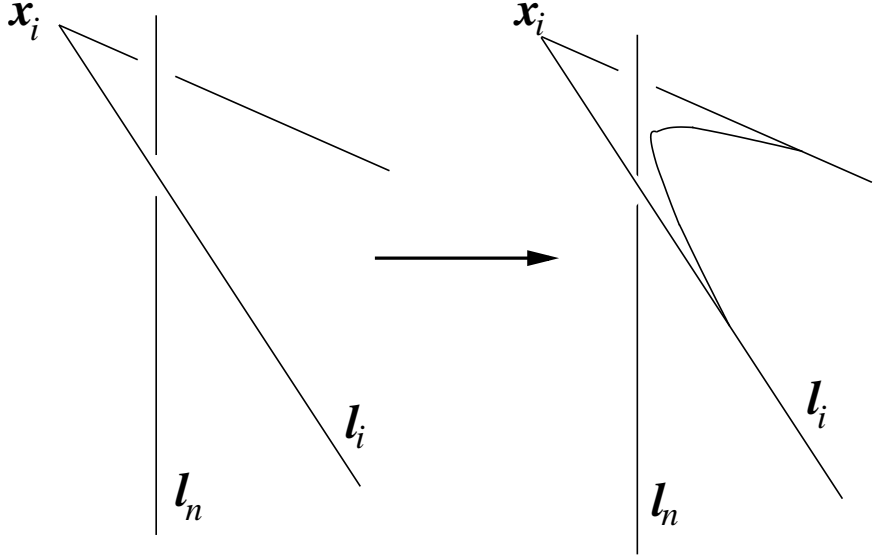}
\caption{A situation that should be avoided: before the smoothing of
  the corner in the point $\bos x_i$, the
  segment $\bos l_n$ was passing under the segment $\bos l_i$. After
  smoothing, the segment $\bos l_i$ has been replaced by an arc of a
  smooth curve in such a way that the segment $\bos l_n$ now passes over
  that arc, potentially changing the topology of the knot.} \label{extrafig} 
\end{figure}

Now we know that the knot obtained after substituting the sharp
corners with smooth arcs of curves is a $G^1-$curve. However, we have to
verify that the topology of the smoothed knot and that of the
original discrete one are the same.
In particular, we have to be sure that, after the replacement of a
corner with an arc of a smooth curve, the dangerous situation depicted
in Fig.~\ref{extrafig} does not occur.
This goal will be achieved by a careful definition of the
lengths of the segments $\bos l_i^\pm$ and $\bos l_i^0$.
To this purpose, we have to derive suitable values $d_i,d'_i$ for 
$i=1,\ldots,N$.
The case $i=N$ is quite an exception, because the segment
$\bos 
l_N$
is followed by $\bos l_1$. This will require a trivial modification
of the procedure that will be presented in the following
and that is valid strictly speaking for $i=1,\ldots,N-1$.
The parameters $d_i,d_i'$ are determined starting from $i=1$
and then proceeding
recursively with the remaining corners in $\bos x_2,\bos x_3,\dots$
At each step $i$, we should check first of all
if $\bos l_i$ and $\bos l_{i+1}$ are
parallel or not. If they are parallel, i.~e. $\frac{\bos l_i\cdot\bos
  l_{i+1}}{l_il_{i+1}}=1$, then no action is required because there is
no sharp corner and it is
possible to pass to the next step $i+1$. In the following, we concentrate
in the treatment of the case in which contiguous segments at the $i-$th joint
are not parallel. 

Starting from $i=1$, the first task consists in finding the 
point $\bos x_{k,1}$ belonging to the knot which is the nearest
to the vertex $\bos x_1$.
The index $k$ refers to the fact that $\bos x_{k,1}$ is lying on a
segment $\bos l_k$ with $k\ne 1,2$. 
The restriction $k\ne 1,2$ is needed to
exclude trivial nearest points belonging to the segments $\bos l_1$ and
$\bos l_2$.

\begin{figure}
\centering
\includegraphics[width=0.9\textwidth]{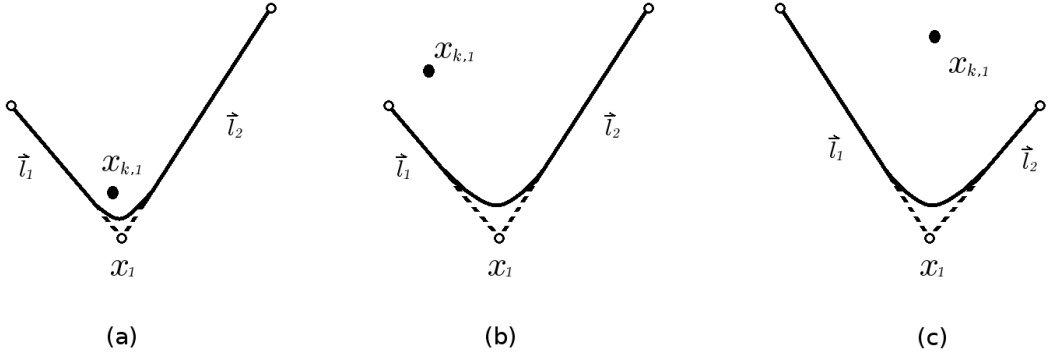}
\caption{This series of pictures illustrates the
meaning of Eq.~(\ref{selone}).
In (a) the nearest point $\bos x_{k,1}$  to
  the vertex $\bos x_1$ is at a distance $d_{k,1}$ from
$\bos x_1$ which is less than
  $\min\left\{\frac{l_1}2,\frac{l_2}2\right\}$.
The length of the segments $\bos l_1^+$ and $\bos l_2^-$ (represented
with dashed lines) is chosen to
be $d_{k,1}$. These segments, which lie inside a sphere $S_{\bos x_1}$ of radius
$d_{k,1}$, are replaced with smooth curves; (b-c)
$\bos x_{k,1}$ is at a distance $d_{k,1}$ from
  $\bos x_1$ which is greater than
$\min\left\{\frac{l_1}2,\frac{l_2}2\right\}$. The length of the
segments $\bos l_1^+$ and $\bos l_2^-$ is chosen to be equal to half
of the length of the shortest among the segments $\bos l_1$ and $\bos
l_2$. Next $\bos l_1^+$ and $\bos l_2^-$ are replaced with smooth curves.
} \label{variouscases} 
\end{figure}

A way in which the position of
$\bos x_{k,1}$ may be computed is presented in the Appendix.
Let us imagine that $\bos x_{k,1}$ is at the distance $d_{k,1}$ from
$\bos x_1$. Then, we choose
\begin{equation}
d_1=d_1'=\min\left\{
d_{k,1},\min\left\{\frac{l_1}2,\frac{l_2}2\right\}
\right\}\label{selone}
\end{equation}
In other words, $d_1$ is set to be equal to $d_1'$. Moreover,
depending on the distance $d_{k,1}$ of the point $\bos x_{k,1}$ from
$\bos x_1$ and on the lengths of the segments $ \bos l_1, \bos l_2$, we
can have the three different possibilities displayed in
Fig.~\ref{variouscases}. 
Fig.~\ref{variouscases}~(b)  refers
to the case in which 
$l_1<l_2$ and
$d_{k,1}\ge\frac{l_1}2$.
Fig.~\ref{variouscases}~(c) shows the analogous situation
in which  $ l_2\le  l_1$ and $d_{k,1}\ge\frac{l_2}2$.
In both cases, by the prescription (\ref{selone}), the
lengths $d_1,d_1'$ can never be greater than half of the length of
the shortest segment between $\bos l_1$ and $\bos l_2$.
When
$d_{k,1}\le \min\left\{\frac{l_1}2,\frac{l_2}2\right\}$, we have
the situation depicted in Fig.~\ref{variouscases}~(a).
Let us notice that the values
$d_1,d_1'$ are selected in such a way that both subsegments $\bos l_1^+$ and $\bos
l_2^-$, together with the smooth arc of curve replacing them, lies
inside a sphere $S_{\bos x_1}$ of radius $d_1$. This sphere contains the point $\bos
x_{k,1}$. Since $\bos
x_{k,1}$ is the point on the knot which is nearest to $\bos x_1$  excluding the
points on the segments $\bos l_1$ and $\bos l_2$, this implies that no
unwanted segment is contained in  $S_{\bos x_1}$. Thus, dangerous
situations
such as those presented in Fig.~\ref{extrafig} are not
possible. 

Now we suppose that all the values of $d_j=d_j'$ have been computed up
to $j<i$. Implicitly, excluding corners in
which the
segments $\bos l_{j-1},\bos l_j$ are parallel, we should assume that the subsegments $\bos
l_j^\pm$ with $j<i$ and $\bos l_i^-$ have already been replaced by the
smooth arcs of 
curves of Eqs.~(\ref{xones}--\ref{xtwos}). 
We also assume that
the smoothing procedure has been carried in such a way that, for $j<i-1$,
the arcs substituting the subsegments $\bos
l_j^+,\bos l_{j+1}^-$ are inside a sphere $S_{\bos x_j}$ of radius $d_j$ and
no other part of the knot after the replacements made so far is
contained in this sphere. The same statement should be true in the
case $j=i-1$ too,
in which the
sphere $S_{i-1}$ of radius $d_{i-1}$ is allowed to contain only the 
arcs of curves which replaced the subsegments
$\bos l_{i-1}^+,\bos l_i^-$.

\begin{figure}
\centering
\includegraphics[width=0.7\textwidth]{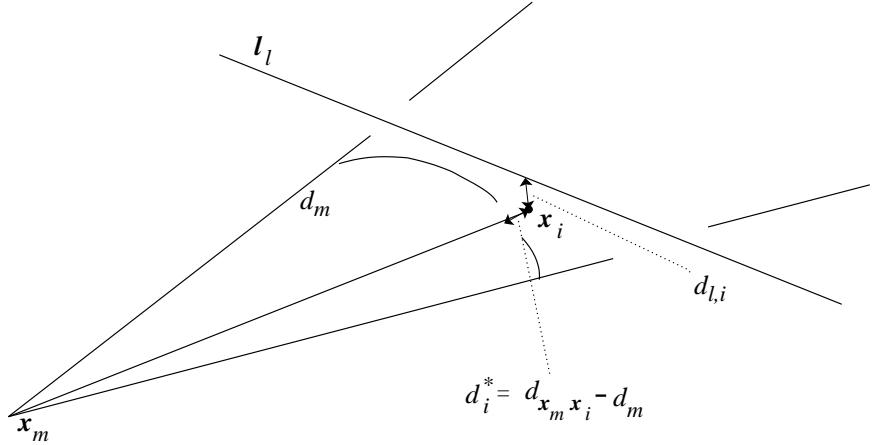}
\caption{This picture explains the meaning of the quantity $d^*_i$
  appearing in Eq.~(\ref{disdef}). It shows also the procedure with
  which the radius $d_i$ of the sphere $S_{\bos x_i}$ is chosen.} \label{problem} 
\end{figure}

At this point we have to
deal with the corner corresponding to the vertex in $\bos x_i$.
As we did for the first corner in $\bos x_1$, we determine the position $\bos x_{l,i}$ of the
point 
which is nearest
to $\bos x_i$
and does not belong to $\bos l_i$ or $\bos l_{i+1}$.
Let's suppose
that $\bos x_{l,i}$ lies on the segment 
$\bos l_l$ with $l\ne i,i+1$ and is at a distance $d_{l,i}$ from $\bos
x_i$. We have also to be sure that no point of the spheres $S_{\bos
  x_k}$ with $k=1,\ldots,i-1$ corresponding
to a corner that has already been substituted is at a distance from
$\bos x_i$  which is smaller than $d_{l,i}$.
To this purpose we compute the minimal distance $d^*_i$ from $\bos
x_i$ to these spheres:
\begin{equation}
d^*_i= \min_{\substack k=1,\ldots,i-1}
\left\{
d_{l,i}, d_{\bos x_k\bos x_i}-d_k
\right\}\label{disdef}
\end{equation}
The meaning of $d^*_i$ is illustrated in Fig.~\ref{problem}. 
 The following three situations should be treated separately:\label{ggg}
\begin{enumerate}
\item $l>i$ and $d^*_i\ge d_{l,i}$. 
\item $l<i$ and $d^*_i\ge d_{l,i}$. 
\item \label{casethree} $d^*_i< d_{l,i}$. 
\end{enumerate}
In the first case, 
the substituted parts of the knot are at a distance $d^*_i$ which is
greater than $d_{l,i}$. As a consequence,
the segment $\bos l_l$ 
contains the point $\bos x_{l,i}$ which is among all points
of the knot the nearest one to $\bos
x_i$. 
Moreover, $\bos l_l$ has not yet been affected by the smoothing
procedure. It is thus possible to
proceed as we did for the first corner $\bos x_1$.
The second case 
is more complicated. The point  $\bos x_{l,i}$ on the segment $\bos
l_l$ is
nearer to $\bos x_i$ than any other point lying on the other segments
or on the parts of the knot that have already been replaced, but the segment $\bos l_l$
has been affected by the smoothing procedure. 
This means that  $\bos x_{l,i}$ could
have been mapped to a new point $\bos
x_{l,i}'$ and it is no longer trivial to determine which is the new nearest point
to $\bos x_i$. 
Three subcases are possible, see Fig.~\ref{safe} for a visual
representation: 
\begin{itemize}
\item[2.-a)] $\bos x_{l,i}$ is near the point $\bos x_{l-1}$ within the
  distance $d_{l-1}$, i.~e.  $\bos x_{l,i}$ lies on the segment $\bos
  l^-_{l}$. This implies that for sure $\bos x_{l,i}$ has been
  already mapped into the point $\bos
x_{l,i}'$ located on the arc of the smooth curve 
that replaced $\bos l_l^-$.
What we know is that both segment $\bos
  l^-_{l}$ and the arc of the smooth curve replacing it are inside the
  sphere $S_{\bos x_{l-1}}$ of radius $d_{l-1}$. The radius $d_i$ of
  the sphere $S_{\bos x_i}$ 
  surrounding $\bos x_i$ should be chosen in such a way that
$S_{\bos x_{l-1}}$ and  $S_{\bos x_i}$ do not penetrate into each
  other. We should also avoid that $S_{\bos x_i}$  penetrates inside
  any other sphere $S_{\bos x_k}$ with $k<i$. Due to the condition
  $d^*_l\ge d_{l,i}$, which is valid in this subcase, this last requirement
  is matched if
   the following inequality is satisfied:
\begin{equation}
d_i,d_i'\le d_{l,i}\label{conddd}
\end{equation}
 Let
  $d_{\bos x_{l-1}\bos x_i}$ be the distance between $\bos x_{l-1}$ and $\bos
  x_i$. Since in our construction $d_{l-1}$ can be only less or equal to the minimal distance between $\bos
  x_{l-1}$ and any segment $\bos l_n$ with $n\ne l-1,l$, it is
  possible to conclude that:
\begin{equation}
d_{\bos x_{l-1}\bos x_i}\ge d_{l-1}\label{eq44}
\end{equation}
If $d_{l-1}$ is strictly smaller than $d_{\bos x_{l-1}\bos x_i}$, then 
the following inequalities hold:
\begin{equation}
0<d_{\bos x_{l-1}\bos x_i}-d_{l-1}<d_{l,i}\label{ineqtobeproved}
\end{equation}
The left inequality is a trivial consequence of our settings.
To prove that $d_{\bos x_{l-1}\bos x_i}-d_{l-1}<d_{l,i}$, we
remember that, under the present 
assumptions:
\begin{equation}
0<d_{\bos x_{l,i}\bos x_{l-1}}<d_{l-1}\label{ineq1}
\end{equation}
 where 
$d_{\bos x_{l,i}\bos x_{l-1}}$ denotes the distance between the points
$\bos x_{l,i}$ and $\bos x_{l-1}$. 
Eq.~(\ref{ineq1}) simply states the fact that the point $\bos x_{l,i}$
lies in the subsegment $\bos l_l^-$ whose length is $d_{l-1}$.
Moreover, using the properties of the
norm expressing the distances on Euclidean spaces, it turns out that 
\begin{equation}
d_{\bos x_{l-1}\bos
  x_i}\le d_{\bos x_{l,i}\bos x_{l-1}}+d_{l,i}\label{ineq2}
\end{equation}
 Applying the second inequality appearing in 
Eq.~(\ref{ineq1}) to (\ref{ineq2}), we obtain  that
$d_{\bos x_{l-1}\bos x_i}-d_{l-1}<d_{l,i}$, thus proving
Eq.~(\ref{ineqtobeproved}).
As a consequence, a sphere $S_{\bos x_i}$ of radius $d_{\bos x_{l-1}\bos x_i}-d_{l-1}$
around the point $\bos x_i$ will never contain any point of the knot
apart from the points belonging to the segments $\bos l_i$ and $\bos
l_{i+1}$.
This is due to Eq.~(\ref{ineqtobeproved}), which states that the
radius
of  $S_{\bos x_i}$ is smaller than the distance $d_{l,i}$ 
between $\bos x_i$ and the
nearest point to $\bos x_i$.
Thus, the condition (\ref{conddd}), which is sufficient to avoid
points of contacts between $S_{\bos x_i}$ and the spheres
$S_{\bos x_k}$ of radii $d_k$ for
$k=1,\ldots,i-1$ and $k\ne l-1$, is fulfilled.
Moreover, when $k=l-1$, by construction
the sphere $S_{\bos x_{l-1}}$ has only one point
of contact with $S_{\bos x_i}$, see Fig.~\ref{safe}.
Thus, it is possible to safely choose $d_i$ and $d_i'$ as follows:
\begin{equation}
d_i=d_i'=\min\left\{
d_{\bos x_{l-1}\bos x_i}-d_{l-1},\min\left\{\frac{l_i}2,\frac{l_{i+1}}2\right\}
\right\}
\end{equation}
With this choice, in fact, no dangerous crossing of lines may occur between
the arcs of curves already replaced and the subsegments that
have still to be treated, see also Fig.~\ref{safe} for a visual
representation of the situation. 
\begin{figure}
\centering
\includegraphics[width=0.7\textwidth]{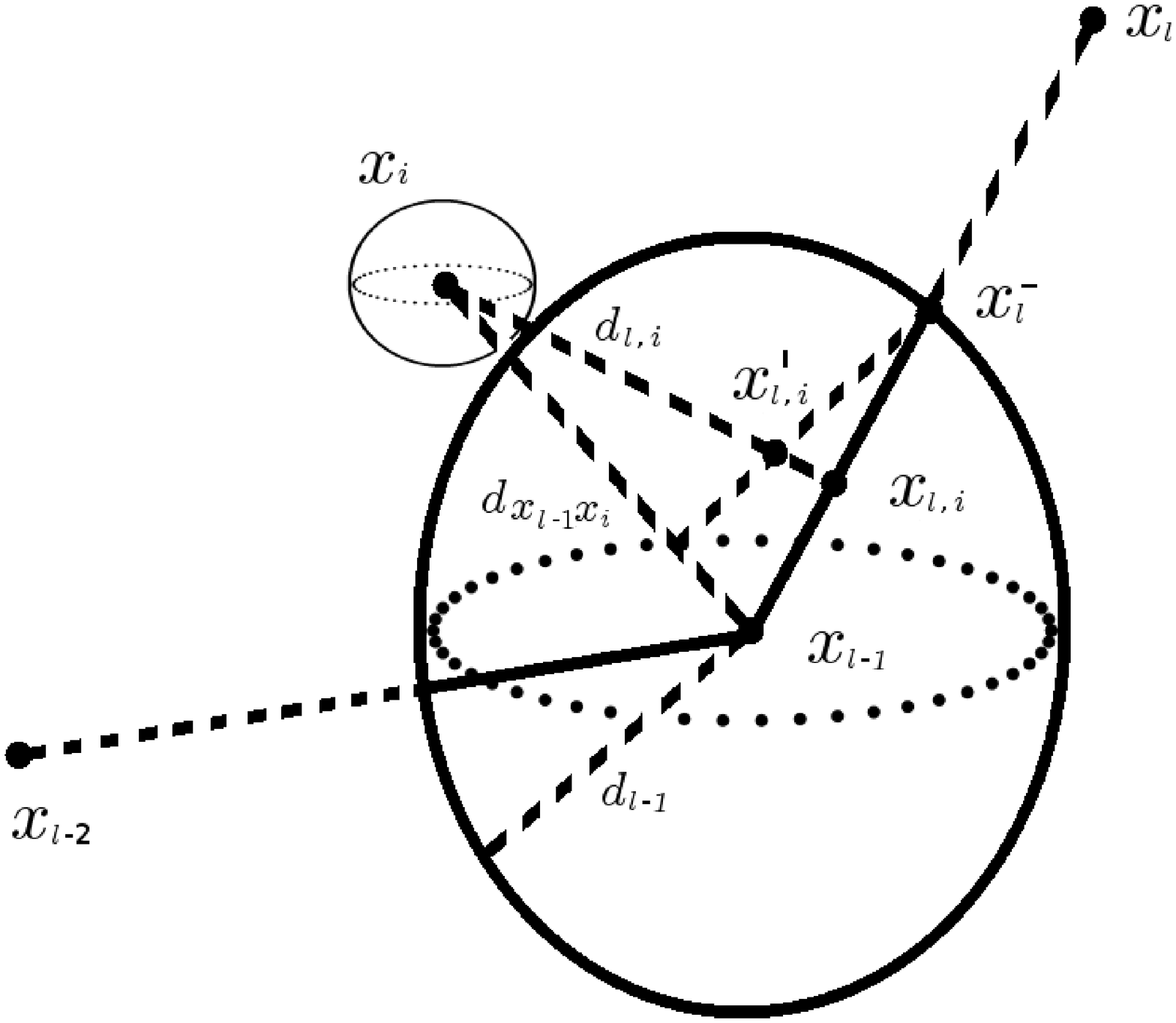}
\caption{This Figure illustrates the subcase 2.-a). We suppose that
  the  nearest point
$\bos x_{l,i}$  to $\bos 
  x_i$ is located on the
  segment $\bos l_l$, with $l<i$. Moreover, $\bos x_{l,i}$ is within
  the distance 
 $d_{l-1}$ from the vertex $\bos x_{l-1}$. Thus, the sharp
  corner at $\bos x_{l-1}$ has been already replaced by the arcs of
smooth curves $\bos
  X^+_{l-1}(S)$ and $\bos X^-_{l}(S)$ of Eqs.~(\ref{xones}) and
  (\ref{xtwos}). The
distance $d_{l-1}$ is smaller by assumption than the distance
  $d_{\bos x_{l-1}\bos x_i}$ from $\bos x_i$ to $\bos x_{l-1}$. 
Moreover, $d_{l,i}$ is smaller than the distance of the point $\bos
x_i$ to the border of any of the spheres $S_{\bos x_k}$ surrounding
the vertices $\bos x_k$ for $k<i$ and $k\ne l-1$.
To be
  safe, the radius of the sphere $S_{\bos x_i}$ around the point $\bos
  x_i$ is chosen 
  to be 
  $d_i=\min\left\{ 
d_{\bos x_{l-1}\bos
  x_i}-d_{l-1},\min\left\{\frac{l_i}2,\frac{l_{i+1}}2\right\} 
\right\}$.} \label{safe} 
\end{figure}

We are left with the particular case $d_{\bos x_{l-1}\bos
  x_i}=d_{l-1}$, which can never happen under the present hypothesis
that $l<i$ and $d^*_i\ge d_{l,i}$. As a matter of fact, if
 $d_{\bos x_{l-1}\bos
  x_i}=d_{l-1}$, 
the quantity $d_i^*$ defined in Eq.~(\ref{disdef}) is zero and, as a
consequence, $d_i^*< d_{l,i}$, because $d_{l,i}$ can never vanish
since the segments $\bos l_l$ and $\bos l_i$ are not allowed to intersect. 
The situations in which $d_i^*=0$ will be treated in the next point
dedicated to case~\ref{casethree} of page~\pageref{ggg}.
\item[2.-b)]  $\bos x_{l,i}$ is near the point $\bos x_l$ within a
  distance which is less than $d_l$, i.~e. it lies on the segment
  $\bos l_l^+$ that has already been substituted.  In that case the
  procedure to determine $d_i,d_i'$ is analogous to that used in the
  case 2.-a).
\item[2.-c)] $\bos x_{l,i}$ is on the subsegment $\bos l_l^0$, In that
  case there is no problem because this subsegment has not been
  replaced with an arc of a curve and we may proceed by taking:
\begin{equation}
d_i=d_i'=\min\left\{
d_{l,i},\min\left\{\frac{l_i}2,\frac{l_{i+1}}2\right\}
\right\}
\end{equation}
\end{itemize}
Finally, we deal with the 
case~\ref{casethree} of page~\pageref{ggg}.
We distinguish the following two subcases:
\begin{itemize}
\item[3-a)] $d_i^*>0$. In this case we choose the radius of the sphere
  $S_{\bos x_i}$ around the point $\bos x_i$ to be equal to $d_i^*$,
  i.~e.:
\begin{equation}
d_i=d_i'=\min\left\{
d_i^*,\min\left\{\frac{l_i}2,\frac{l_{i+1}}2\right\}
\right\}
\end{equation}
\item[3-b)] $d_i^*=0$.  
As a consequence, there exist a sphere $S_{\bos x_m}$ of
radius $d_m$ centered at the point $\bos x_m$ for some value of
$m\le i-1$  such that
 $\bos x_i$ lies on the border of this sphere.
On the
  other side, no point of the knot may be nearer to $\bos x_i$ than
  $\bos x_{l,i}$, which is at the distance $d_{l,i}>d^*_i=0$ from $\bos
  x_i$.
Yet, it may happen that there exists another sphere, let say $S_{\bos
  x_n}$, with $n\le i-1$ and $n\ne m$ centered at a point $\bos x_n$,
such that  the distance between the border of $S_{\bos
  x_n}$ and $\bos x_i$ is less than $d_{l,i}$.
To take into account this case, it is useful to
compute also the quantity $d_i^{*\prime}$\footnote{We do not discuss
  here the particular case in which $d_i^{*\prime}$ is also equal to
  zero, which can be easily treated.}:
\begin{equation}
d^{*\prime}_i=\min_{\substack {k=1,\ldots,i-1\\k\ne m}}\left\{
d_{l,i}, d_{\bos x_k\bos x_i}-d_k
\right\}\label{disdeftwo}
\end{equation}

At this point a possible way to define $d_i$ and $d_i'$ consists in putting
\begin{equation}
d_i=d_i'=\min\left\{
\min\left\{d_i^{*\prime},\frac{d_k}2\right\},\min\left\{\frac{l_i}2,\frac{l_{i+1}}2\right\}\label{ecns}
\right\}
\end{equation}
and to decrease the value of the radius of the sphere $S_{\bos x_k}$
as follows:
\begin{equation}
d_{k}\longrightarrow \frac{d_{k}}2\label{ecdlmo}
\end{equation}
Clearly, the choice of $d_i,d_i'$
given above does not allow crossings of the lines that can change the
type of the knot. 
\end{itemize}
Let us stress that, within the procedure 
illustrated above, at each stage $i$ the
arcs of the curves replacing the corners
are always contained inside a sphere centered at the
$i-$th corner. The spheres corresponding to different corners never
intersect themselves apart from one point on their surfaces which may be in
common. In this way, there is no possibility that the trajectory of
the knot crosses itself during the replacement of the sharp corners
causing unwanted changes of the topology of the knot.
The result is that the trajectory of the discrete
knot becomes a $G^1-$curve, which is the condition sufficient in order
to compute the Vassiliev knot invariant of degree 2 without the
systematic errors related to the presence of the sharp corners.
The only drawback is that the parameters $d_i$, and thus the portions
of the corners that have been replaced, become smaller and
smaller when the value of $i$ increases. This is a drawback because  an
extensive Monte 
Carlo sampling procedure  is required
in order to evaluate $\varrho(C)$ with a good approximation if the
arcs of
curves replacing the corners exhibit sharp turns.

The method explained above is somewhat complicated, but this is only
because we are considering 
very general discrete knots, 
defined off or on~lattice and
with segment of different lengths. On a lattice, many simplifications
are possible. For instance, 
on a simple cubic lattice one may always choose $d_i=d_i'=\frac{l_i}2$
without the risk of changing the topology of
the knot. In that case it is easy to check that
Eqs.~(\ref{xones}--\ref{xtwos}) reduce to the equation:
\begin{equation}
\bos X_i^\pm(S)=\left(1-\sin\left(\theta_0^\pm(S)\right)\right)\bos l_i^++
\left(1-\cos\left(\theta_0^\pm(S) \right)\right)\bos l_{i+1}^-+\bos x_i
\end{equation}
where
\begin{equation}
\theta_0^+(S)=\left(S-[S]-\frac12\right)\frac\pi2
\end{equation}
and
\begin{equation}
\theta_0^-(S)=\left(S-[S]+\frac12\right)\frac\pi2
\end{equation}
Similar simplifications occur in the ansatz of
Eqs.~(\ref{xones}--\ref{xtwos}) in 
the particular case (off and on 
lattice) in which all the segments have the same length and the
$d_i$'s coincide with the $d_i'$'s.

\begin{figure}
\centering
\includegraphics[width=0.7\textwidth]{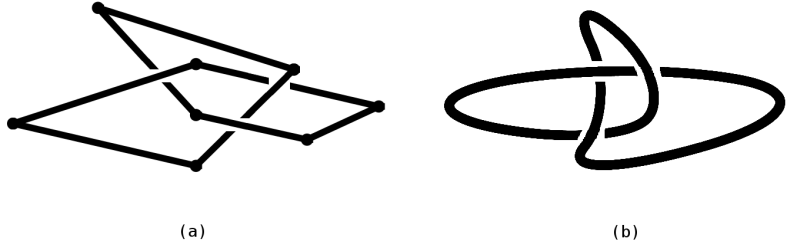}
\caption{  A knot $3_1$ with minimal length defined off lattice
before and after the smoothing procedure.} \label{smoothing}
\end{figure}

Finally, after the values $d_i,d_i'$ are computed for all the corners, the
points on the knot are sampled using the following prescriptions: 
\begin{description}
\item[Step 1] Pick up a random number $S$ in the interval
  $[0,N]$. These random variables are necessary in order to generate
  the points $\bos X(S)$ needed for the Monte Carlo procedure. The
  value of $S$ identifies a segment $\bos l_i$ with $i=[S]+1$ and end
  points $\bos x_i=\bos 
  X([S]+1)$ and $\bos x_{i-1}=\bos
  X([S])$. This procedure works also for $i=1$ provided the point $\bos x_0$ is identified with $\bos x_N$.
\item[Step 2] 
We assume that the curve is oriented  in such a way that  the $i-$th
segment $\bos l_i$ is coming before the segment $\bos l_{i+1}$. 
Now we have to check
if one of the following three conditions are satisfied: 
\begin{eqnarray}
0 \le &S-[S]& < \frac{d_{i-1}'}{ l_i} \label{condi1}\\
\frac{d_{i-1}'}{ l_i} \le &S-[S]& < \frac{ l_i-d_i}{ l_i} \label{condi2}\\
\frac{ l_i-d_i}{ l_i} \le &S-[S]& < 1 \label{condi3}
\end{eqnarray}
The first condition (\ref{condi1}) identifies the subsegment $\bos
l_i^-$, the second condition (\ref{condi2}) the subsegment $\bos
l_i^0$ and the third one (\ref{condi3}) the subsegment $\bos l_i^+$. 
\item[Step 3] When condition (\ref{condi1}) is fulfilled,  
verify if the relation
\begin{equation}
\frac{\bos l_{i-1}^+\cdot{\bos l_i^-}}
  { l_{i-1}^+ { l_i^-}}=-1
\end{equation}
is satisfied.
If yes, the segments $\bos l_{i-1}^+$ and $\bos l_i^-$ around the
corner $\bos x_{i-1}$ 
are antiparallel and  
there is no sharp corner to be smoothed up. In that
case the smoothed curve coincides with the old one and the
parametrization given in 
Eqs.~(\ref{Xsdef}--\ref{Xsdefadd}) is still valid. 
If instead the segments $\bos l_{i-1}^+$ and $\bos l_i^-$ are not
antiparallel, then we have to use the prescription
of Eq.~(\ref{xtwos}), which is valid on $\bos l_i^-$ after the
replacement $i\longrightarrow i-1$.
The point $\bos X(S)$ is projected onto the point $\bos X_i^-(S)$ with
the help of Eq.~(\ref{xtwos}).
A similar procedure is adopted when condition
(\ref{condi3}) is true. In that case the condition of being
antiparallel is concerning the segments $\bos l_i^+$ and $\bos l_{i+1}^-$: 
\begin{equation}
\frac{\bos l_i^+\cdot{\bos l_{i+1}^-}}
  {l_i^+ {l_{i+1}^-}}=-1
\end{equation}
If it turns out that $\bos l_i^+$ and $\bos l_{i+1}^-$ are not
antiparallel, then the point $\bos X(S)$ should be mapped 
using the curve in Eq.~(\ref{xones}).
Finally, when condition (\ref{condi2}) is satisfied, we are on the
subsegment $\bos l_i^0$ away from any corner. As a consequence, for
the values of $S$ in the interval $[\frac{d_{i-1}'}{ l_i},\frac{
    l_i-d_i}{ l_i}]$, it is possible to apply the old
parametrization of 
Eqs.~(\ref{Xsdef}--\ref{Xsdefadd}).
\end{description}
An example of curve describing a discrete knot $3_1$ off lattice  
before and after the smoothing procedure is shown
in Fig.~\ref{smoothing}.

The smoothing procedure illustrated above has been applied to many
examples of different knots of various lengths. It delivers values of
the Vassiliev invariant of degree 2
which are approaching the exact value with a precision that increases with
the increase of the number of samples used in the Monte Carlo
integration algorithm.
In Table~\ref{table2} we report for instance the case of a knot $3_1$ with
  $24$ segments computed using gradually increasing numbers of
samples.

\begin{table}
\centering
\begin{tabular}{|c| c|}
\hline
 $n$ & $\varrho(C)$ \\ [0.5ex]
\hline
$10^6$ &$1.9096\pm 0.0991$ \\
$10^7$ &$1.9179\pm 0.0326$ \\
$10^8$ &$1.9170\pm 0.0095$ \\ 
$10^9$ &$1.9168\pm 0.0032$ \\ 
 [1ex]
\hline
\end{tabular}
\caption{Computation of the knot invariant $\varrho(C)$ for the knot $3_1$ with
  $24$ segments on a simple cubic lattice. The results of the
numerical calculation of $\varrho(C)$
are displayed for different values
 of
 the number of samples $n$ used in the Monte Carlo integral
  procedure. As it can be seen, 
by gradually increasing $n$,
the numerical values of $\varrho(C)$
asymptotically approach the analytical value of the
Vassiliev knot invariant of degree 2 which,
in the case of the knot $3_1$,
 is approximately equal to $1.9167$.}  
\label{table2}
\end{table} 
Table~\ref{table1} illustrates how the presence of the sharp corners
affects
the calculations of $\varrho(C)$.
\begin{table}[ht]
\centering
\begin{tabular}{|c| c| c| c| c|}
\hline
 knot type & $\varrho_a(C)$ & $\varrho_{sp}(C)$ & $\varrho_{ns}(C)$ &$n_{sc}$\\ [0.5ex]
\hline
$0_1$ &$-\frac{1}{12}$   &$-0.0839\pm 0.0332$ &$+0.5526\pm 0.0569$ &77\\
$3_1$ &$+\frac{23}{12}$  &$+1.9170\pm 0.0553$ &$+2.4781\pm 0.0465$ &68\\
$4_1$ &$-\frac{25}{12}$  &$-2.0847\pm 0.0533$ &$-1.5214\pm 0.0845$ &68\\ 
$5_1$ &$+\frac{71}{12}$  &$+5.9174\pm 0.0653$ &$+6.4523\pm 0.0845$ &65\\ 
$6_1$ &$-\frac{49}{12}$  &$-4.0856\pm 0.0723$ &$-3.5717\pm 0.1007$ &62\\ 
$7_1$ &$+\frac{143}{12}$ &$+11.9173\pm 0.0652$ &$+12.4258\pm 0.1217$ &62\\ 
$8_1$ &$-\frac{73}{12}$  &$-6.0822\pm 0.0529$ &$-5.6380\pm 0.0774$ &54\\ 
$9_1$ &$+\frac{239}{12}$ &$+19.9158\pm 0.0855$ &$+20.4041\pm 0.1579$ &59\\ [1ex]
\hline
\end{tabular}
\caption{ This table provides the values of
the Vassiliev knot invariant of degree 2
  for the knots $0_1$, $3_1$, $4_1$, $5_1$, $6_1$, $7_1$, $8_1$
  and $9_1$. $\varrho_a(C)$ denotes the analytical value of the
  knot invariant.  $\varrho_{sp}(C)$ refers to the
   results of the computation of the
 knot invariant obtained after performing the smoothing procedure 
 described in Section \ref{section4}.
  $\varrho_{ns}(C)$ is instead the value of the knot invariant
  derived without the smoothing
   procedure. The data of $\varrho_{sp}(C)$  
and  $\varrho_{ns}(C)$ have been
computed using the same number of samples, which varies depending on
the kind of knot. Finally, $n_{sc}$ is the number of sharp corners contained in
the knot before the smoothing procedure.}  
\label{table1}
\end{table} 
In the second column of Table~\ref{table1}, the outcome
$\varrho_a(C)$ of the  analytical
computation of $\varrho(C)$ is provided
for several knots with number of segments $N=90$. Within the given
errors, the values $\varrho_{sp}(C)$
obtained by Monte Carlo integration after the smoothing procedure
(sp)
 are in
agreement with the analytical results, see the third column of
Table~\ref{table1}.
We are also reporting the upshot of the calculations performed
without the smoothing procedure, see the values of $\varrho_{ns}(C)$ in the
fourth column of
Table~\ref{table1}. 
The differences between $\varrho_{sp}(C)$  and $\varrho_{ns}(C)$ show
that indeed the presence of sharp corners in the case of discrete
knots does not allow the correct evaluation of the knot invariant
$\varrho(C)$.

\section{Speeding up the Monte Carlo algorithm} \label{section5}
The computation of the Vassiliev invariant of
degree 2  by Monte Carlo methods is much more convenient
than by traditional numerical techniques.
For instance, in order to evaluate $\varrho(C)$ with sufficient
precision
in the case of knots of length $L\le 120$, a few millions of
samples are enough. This is a quite good performance if we take into
account that, for $N=120$ the total volume to be checked is 
$\frac{120^4}{24}\sim 9\cdot 10^6$.
If a very high precision is required or $N$ is large, 
the sampling procedure can be easily parallelized on a computer.
Still, the numerical evaluation of $\varrho(C)$ becomes challenging
especially in the case of  knots consisting of
a large number of segments and it is advisable to adopt some strategy in
order to reduce the calculation time.
Let us notice at this point that,
in practical applications, knot invariants are mainly used
in order to make assessments on the topological configuration of a
knot which is a priori unknown.
To that purpose, it is not necessary to evaluate
$\varrho(C)$ beyond a certain precision, as it will be evident
from the following discussion.
First of all, let's recall the fact that
there is no knot invariant that is able to distinguish
unambiguously all different types of knots.
The  Vassiliev invariant of degree 2 is not an exception to this
rule, but still may be considered as a relatively powerful knot
invariant. For example, it is able to distinguish uniquely
the knots $9_1$ and $10_3$ from all other knots up to
ten crossings. Of course, there are many knots for which
the second coefficient of the Conway polynomial
$a_2(C)$ is the same. 
This implies that $\varrho(C)$, which is related to
$a_2(C)$ by Eq.~(\ref{rela2rho}),
can at most be used to distinguish
classes of knots having different values of  $a_2(C)$.
A nice characteristic of $\varrho(C)$ consists in the fact that, if
two knots $C$ and $C'$ can be resolved by it, then the condition
\begin{equation}
\left|
\varrho(C)-\varrho(C')
\right|\ge 2\label{ooo}
\end{equation}
is always satisfied.
As a consequence, in order to ascertain the difference between two  knots
with the help of a Monte Carlo calculation of $\varrho(C)$, it is not
necessary to push the standard deviation  $\sigma$ of the numerical
calculation of $\varrho(C)$ 
below a given threshold value $\sigma_{threshold}$. For instance, we
can choose:
\begin{equation}
\sigma_{threshold}=\frac 1{6.11}\sim 0.16\label{varcond}
\end{equation}
If $\sigma=\sigma_{threshold}$, in fact, the probability that the
Monte Carlo evaluation of $\varrho(C)$ gives a result
within an error of $\pm 1$ or greater is of the order $1\cdot 10^{-9}$,
i.~e. this event is very unlikely. If two knots $C$ and $C'$ are
distinguishable by using the Vassiliev invariant of degree 2,
then, due to  Eq.~(\ref{varcond}),
this precision is more than enough to state if  $C$ and $C'$ are
different or not with a satisfactory confidence level.
The possibility of putting a lower cutoff to
the standard deviation as in Eq.~(\ref{varcond}) is very helpful in
practical calculations because,
in order to decrease the standard deviation, it is necessary
to increase the number of samples used in the Monte Carlo integration
method. Of course, this leads  to a consistent increase of the
calculation time. In fact, the time $\tau$ necessary for computing
$\varrho(C)$ scales linearly with the number of samples $n$, but an
increasing of $n$ by a factor $\lambda>1$ produces an improvement of
$\sigma$ only by a factor $\frac1{\sqrt{\lambda}}$,
i.~e. $\sigma\longrightarrow \frac\sigma{\sqrt{\lambda}}$. We have
checked that this scaling law, that is predicted in the case of
gaussian distributions, is actually verified in the present context.

\begin{figure}
\centering
\includegraphics[width=0.5\textwidth]{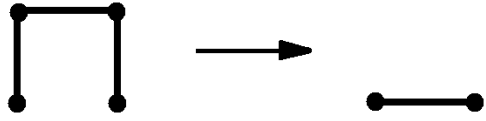}
\caption{Reduction to a single segment of an element of the discrete
  knot composed by
three consecutive segments such that its
   ends are at a distance equal to the 
  unit size on a simple cubic lattice. The topology of the knot is not
  affected by this 
reduction.} \label{tadpoles} 
\end{figure}

\begin{figure}
\centering
\includegraphics[width=0.5\textwidth]{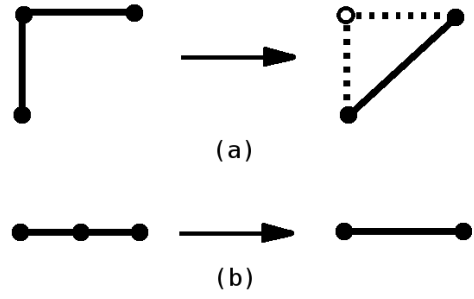}
\caption{This figure displays the possible configurations on a simple
  cubic lattice 
of two contiguous segments and their substitution with a single
segment. The topology of the knot is left unchanged after the
substitution.} \label{two} 
\end{figure}

\begin{figure}
\centering
\includegraphics[width=0.5\textwidth]{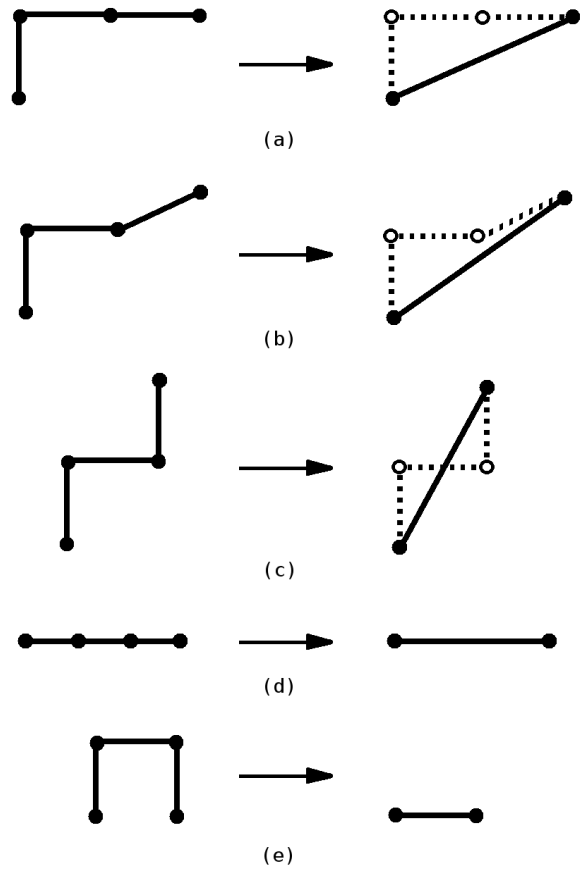}
\caption{This figure displays the possible configurations on a cubic lattice
of three contiguous segments and their substitution with a single
segment.} \label{three}
\end{figure}

Besides the standard deviation, another important factor which
determines the computation time $\tau$ is the number $N$ of segments
(or equivalently the number of arcs of $G^1-$trajectories) in which
the trajectory of the knot is realized. Due to the presence of a
quadruple contour integral in $\varrho_2(C)$, see
Eqs.~(\ref{rhotwomi}) and (\ref{F2def}), $\tau$ scales with respect to
$N$ according to the fourth power, i.~e. $\tau\propto N^4$.
We stress the fact that $\tau$ depends on the number of points $N$ and
not
on the knot length  $L$. 
The length of the knot depends in fact on the
lengths $l_i$
 of the segments $\boldsymbol l_i$ for  $i=1,\ldots,N$. In turn, the $l_i$'s
can be made 
arbitrarily small by the rescalings  
\begin{equation}
l_i\longrightarrow
l_i'=\eta l_i \qquad\qquad
i=1,\ldots,N \label{rescaling}
\end{equation}
where $\eta$ is any real parameter
such that $0<\eta<1$. In fact, rescalings of this kind do not affect
the value of $\varrho(C)$, which is a scale invariant quantity.
To prove the statement that $\tau\propto N^4$, it is sufficient to
decompose, as it 
has been done in Ref.~\cite{YZFFJSTAT2013},
the quadruple integral in
 Eq.~(\ref{rhotwomi}) into a quadruple sum 
of contour integrals in which  the contours are
the segments $\boldsymbol l_i$ themselves. 
With the rescalings  of Eq.~(\ref{rescaling}) it is then possible
to reduce the lengths of these
segments to infinitesimal quantities. As a consequence,
the contour integrals
over the  segments $\boldsymbol l_i$ can be computed exactly.
In this way, only  the quadruple sum over the indexes labeling the
segments remains which contains exactly $\frac {N^4}{24}$ terms.
For knots composed by a large number of segments, 
this implies that
the number of samples necessary for obtaining a
satisfactory result from the Monte Carlo integration algorithm becomes
prohibitively high. This problem can be partially avoided by adopting
procedures that are able to shorten a given discrete knot reducing its number
of segments without changing its topology.
In the following, a few such procedures will be proposed, most of them
valid on a simple 
cubic lattice before applying the smoothing procedure:
\begin{enumerate}
\item
For a general discrete polymer, it is always possible to group
together contiguous segments that are parallel.
\item On a simple cubic lattice, configurations of three segments
  whose ends are at a distance equal to the size of a lattice edge can
  be easily substituted by one segment as shown in
  Fig.~\ref{tadpoles}.
This reduces the length of the knot by two segments every time this
configuration is encountered.
\item
Always on a simple cubic lattice it is possible to group together two
or three contiguous segments in a single one, see Figs.~\ref{two} and
\ref{three}. 
We note that the first substitution in Fig.~\ref{three} can cause
intersections between two segments after the grouping and should be
treated with some care.
\end{enumerate}
Other algorithms to decrease the size of a knot can be found in
Refs.~\cite{KoniarisMuthukumar} and \cite{TaylorAsrodi}, where the KMT
radiation scheme has been introduced.

Finally, we show yet another strategy which can be used to
speed up the Monte Carlo procedure when the knot invariant is needed
to assess the topology of a knot after a random transformation that can
be in principle not topology-preserving. 
We will assume that the random transformation involves $K$ contiguous
segments, where $0\le K\le N$.
Instead of calculating the whole knot
invariant $\varrho(C)$, it is much better to compute the difference
$\Delta\varrho(C)=\varrho(C_T)-\varrho(C_R)$. Here $C_R$ is
the polymer conformation before the transformation. By hypothesis, 
$C_R$ is in the desired topological configuration.
$C_T$ denotes instead the conformation after the transformation.
The problem is to ascertain if the new knot $C_T$ is topologically
equivalent to $C_R$. 
Clearly, if $\varrho(C)$ would be a perfect topological invariant,
this would be true only if $\varrho(C_T)-\varrho(C_R)=0$.
Unfortunately, $\varrho(C)$ is not able to distinguish unambiguously 
two different topological configurations. However, 
$\varrho(C)$ is a quite powerful knot invariant. Moreover, the
probability to pass with a single random transformation from a knot
$C_R$ to a knot $C_T$ such that $C_R$ and $C_T$ are topologically
inequivalent but are still characterized by the same value of the
Vassiliev knot invariant of degree 2 seems not to be quite high.
For this reason it is licit to expect that $C_T$ and $C_R$ are very
likely to be of the same topological type if the result of the
numerical evaluation of the difference
$\Delta\varrho(C)=\varrho(C_T)-\varrho(C_R)$ 
gives a nearly vanishing result.
The advantage of considering 
$\Delta\varrho(C)$ 
instead of computing the whole value of the knot invariant for $C_T$ is that 
in this way we can ignore the part of the knot that has remained
unchanged after the transformation. 

To realize that this strategy is convenient, we limit ourselves to the
calculation of $\varrho_2(C)$, which 
is the contribution to $\varrho(C)$ 
that requires the biggest computational effort. From
Eq.~(\ref{vol2}), it turns out that to estimate the value of
$\varrho_2(C)$,  the volume to be explored via the Monte Carlo
sampling procedure is equal to
 $N^4/24$.  If the number of the changed segments is $K$,
then at the leading order the number of terms involving only segments
that have  not been affected by the transformation is equal to
$(N-K)^4/24$.  
As a consequence, the number of summands to be taken into account in
the evaluation of the difference $\varrho(C_T)-\varrho(C_R)$ is  
\begin{eqnarray}
S_K&=&\frac{N^4}{24}-\frac{(N-K)^4}{24}\nonumber\\&=&\frac{NK^3}{6}-
\frac{N^2K^2}{4}+
\frac{N^3K}{6}-\frac{K^4}{24}  
\label{snumber}
\end{eqnarray}

Clearly, the minimum of $S_K$ with respect to $K$ occurs when $K=4$ (we
do not consider here  transformations with less than $4$
segments). Due to the fact that the derivative of $S_K$ with respect
to $K$ in the 
range $0\le K\le N$ is always positive, because
$\frac{dS_K}{dK}=\frac{(N-K)^3}{6}>0$, it turns out that $S_K$ grows
with $K$ 
until it reaches its maximum when $K=N$ (transformations of more
than $N$ segments do not make sense). If $K$ is small with respect to
$N$, Eq.~(\ref{snumber})
shows that the volume to be explored in the computation of
$\Delta\rho_2(C)$ is much less than that needed to
obtain $\varrho_2(C)$. For instance, when $K=N/5$, we obtain: 
\begin{eqnarray}
S_K=0.5904\times\frac{N^4}{24}
\end{eqnarray}
so that only $60\%$ of the original volume $\frac{N^4}{24}$ should be
considered. In the best case, $K=4$, instead,
\begin{eqnarray}
S_4\sim\frac{2N^3}{3}\label{53}
\end{eqnarray}
which implies an enormous gain in speed.

\section{Conclusions}\label{section6}

In this work an algorithm to compute the Vassiliev invariant of degree
2 $\varrho(C)$ for any knot $C$ has been provided.  Particular
attention has been devoted to discrete knots whose trajectories
consist of segments connected together at their ends forming
$C^0-$curves in the space. This case is
relevant in numerical simulations, where knots are
forcefully discretized and represented as $C^0-$curves.
The Vassiliev invariant of degree 2 is probably the simplest knot
invariant that can be defined in terms of multiple integrals computed
along the contour of the knot itself. For this reason, there are
chances that in the future this invariant will play for knots the same
role played by the Gauss linking number in numerical studies of links
formed by ring-shaped quasi one-dimensional objects. A suitable
parametrization of  discrete  knots has been introduced, see
Eqs.~(\ref{equ11}--\ref{Xsdefadd}) and the problem of computing the
multiple contour integrals has been tackled using the Monte
Carlo integration scheme summarized by the general equation
(\ref{mcf}). 
In principle, using the scale invariance of $\varrho(C)$,
it is possible to reduce arbitrarily the length of the knot,
because a change of scale does not alter its topology.
However, even in the best case in which the lengths of the segments
become infinitesimal, so that the integration over them can be easily
approximated, standard methods require the computation of a
sum over $\frac{N^4}{24}+\frac{N^3}{6}$ terms in order to obtain the
value of  $\varrho(C)$, as shown in \cite{YZFFJSTAT2013} and in the
previous Section.
With the Monte Carlo integration a considerably smaller number of
samples is necessary in order to evaluate $\varrho(C)$ with a
satisfactory approximation.
For this reason, the standard integration methods like the Simpson's rule
are decidedly slower as it is pointed out by
 the example provided at the beginning of Section
\ref{section5} of a knot with $120$ segments.
Further refinements of the naive Monte Carlo integration scheme
presented in Section \ref{section3} are difficult to be implemented or do not
increase substantially the speed of the computations. For instance,
the division of the integration domain, which is recommended as 
one of the strategies
 to improve the sampling efficiency, 
did not  lead to significant improvements.
On the other side, it is not easy to guess which distribution of
the values of the integration variables could be suitable in order
to enhance the sampling procedure.

In the computations of  $\varrho(C)$ for discrete
knots consisting of a set of segments, we have found that the results
never coincide with the analytical values. This is
expected because of the sharp corners at the points where the segments
join together. There is some correlation between the
number of these corners and the systematic error in the evaluation
of $\varrho(C)$ that apparently
depends on
the type of the knot and the number of its segments $N$. However, to
establish a general relation between that systematic error
and $N$ which could be valid for every knot has not been possible.
To solve the problem of the sharp corners, a procedure for the smoothing
of discrete knots has been presented in Section \ref{section4}. This
procedure 
transforms the discrete knot into a $G^1-$curve. After the smoothing,
it has been possible to evaluate $\varrho(C)$ with an arbitrary
precision by gradually increasing the number of samples used in the
Monte Carlo integration of the multiple integrals entering the
definition of $\varrho(C)$, see for example
Table~\ref{table2}. Despite the advantages of the adopted 
Monte Carlo method with respect to the traditional integration
techniques, the calculation of $\varrho(C)$ becomes challenging when
the number of segments composing the knot is large.
For certain physical applications of the knot invariant $\varrho(C)$,
we show that the time necessary for its evaluation, which
approximately scales as the fourth power of $N$ when $N$ is large, can
be reduced in such a way that scales with the third power of $N$, see
Eq.~(\ref{53}). Moreover, we present an algorithm to reduce the number
of segments of a knot defined on a cubic lattice by a factor three
without changing the topology. After this reduction, the knot is no
longer defined on a lattice, but still the general smoothing procedure
of Section \ref{section4} and the provided Monte Carlo integration scheme of
Section \ref{section3} can be applied in order to obtain the value of
$\varrho(C)$.
To give an idea of the efficiency of the methods for reducing the
number of segments explained in Section \ref{section5}, in the case of
a knot with $N=1000$ originally constructed on a simple cubic lattice,
the number of segments in the final configuration obtained after the
treatment ranges between 255 and 300 depending on the initial shape of
the knot.
Finally, it is important to notice that, in order to distinguish the
topology of different knots, it is not necessary to achieve a standard
deviation that is lower than the threshold value given in
Eq.~(\ref{varcond}). In this case, in fact, the probability that one
knot could be confused with a topologically inequivalent one due to
statistical errors is very small, of the order of $10^{-9}$.
Work is in progress in order to generalize the methods presented in
this work to the case of the triple invariant of Milnor that describes
the links formed by three knots.
\acknowledgments { The support of the Polish National Center of Science,
scientific project No. N~N202~326240, is gratefully acknowledged.
The simulations reported in this work were performed in part using the HPC
cluster HAL9000 of the Computing Centre of the Faculty of Mathematics
and Physics at the University of Szczecin. }

\appendix
\section{}
To compute the nearest point $\bos x_{k+1,i}$ of a segment $\bos l_{k+1}=\bos x_{k+1}-\bos x_k$ from the vertex $\bos x_i$, we pick up on $\bos l_{k+1}$ a general point $\bos X_{k+1}(\sigma)$ as follows:
\begin{equation}
\bos X_{k+1}(\sigma)=\bos x_k +(\bos x_{k+1}-\bos x_k)\sigma\label{d1}
\end{equation}
with $\sigma\in [0,1]$. The distance between this point and $\bos x_i$ is $\sqrt{(\bos X_{k+1}(\sigma)-\bos x_i)^2}$.
If $\bos X_{k+1}(\sigma)$ is the nearest point to $\bos x_i$, then it satisfies the condition
\begin{equation}
\frac{d \sqrt{(\bos X_{k+1}(\sigma)-\bos x_i)^2}}{d \sigma}=0\label{d2}
\end{equation}
Inserting Eq.~(\ref{d1}) in (\ref{d2}) and solving Eq.~(\ref{d2}) with respect to $\sigma$, we obtain that the point of $\bos l_{k+1}$ at the minimal distance from $\bos x_i$ corresponds to the following value of $\sigma$:
\begin{equation}
\sigma_{\min}=-\frac{(\bos x_{k+1}-\bos x_k)\cdot (\bos x_k-\bos x_i)}{(\bos x_{k+1}-\bos x_k)^2}
\end{equation}
Three cases may occur:
\begin{itemize}
\item[1)] If $\sigma_{\min} \ge 1$, this means that the nearest point occurs on 
a line having the same direction of $\bos l_{k+1}$ at a distance $\sigma_{\min} \ge 1$ from the point $\bos x_k$. This means that the nearest point to $\bos x_i$ on the segment $\bos l_{k+1}$ is $\bos x_{k+1,i}=\bos x_{k+1}$ and its distance from $\bos x_i$ is $d_{k+1,i}=|\bos x_{k+1}-\bos x_i|$.
\item[2)] If $\sigma_{\min}\le 0$, the nearest point occurs on a line having the same direction of $\bos l_{k+1}$ at a distance $-\sigma_{\min}$ from $\bos x_k$. On the segment $\bos l_{k+1}$, the nearest point is in this case the point $\bos x_k$. Its distance from $\bos x_i$ is $d_{k+1,i}=|\bos x_k-\bos x_i|$.
\item[3)] If $0<\sigma_{\min}<1$, then $\bos x_{k+1,i}$ lies on the segment $\bos l_{k+1}$ and $\bos x_{k+1,i}=\bos x_{k+1} +(\bos x_{k+1}-\bos x_k)\sigma_{\min}$. The distance of $\bos x_{k+1,i}$ from $\bos x_i$ is in this case:
\begin{equation}
d_{k+1,i}=\sqrt{(\bos x_k-\bos x_i)^2-\frac{\left[ (\bos x_{k+1}-\bos x_k)\cdot (\bos x_k-\bos x_i)\right]^2}{(\bos x_{k+1}-\bos x_k)^2}}
\end{equation} 
\end{itemize}
By repeating this procedure for all segments $\bos l_{k+1}$ with $k=0,\cdots,N-1$ and $k\ne i-1, i$, we obtain the location of the point of the knot which is not belonging to $\bos l_i$ and $\bos l_{i+1}$, and it is the nearest one from $\bos x_i$.


\begin{thebibliography}{3}%
\makeatletter
\providecommand \@ifxundefined [1]{%
 \@ifx{#1\undefined}
}%
\providecommand \@ifnum [1]{%
 \ifnum #1\expandafter \@firstoftwo
 \else \expandafter \@secondoftwo
 \fi
}%
\providecommand \@ifx [1]{%
 \ifx #1\expandafter \@firstoftwo
 \else \expandafter \@secondoftwo
 \fi
}%
\providecommand \natexlab [1]{#1}%
\providecommand \enquote  [1]{``#1''}%
\providecommand \bibnamefont  [1]{#1}%
\providecommand \bibfnamefont [1]{#1}%
\providecommand \citenamefont [1]{#1}%
\providecommand \href@noop [0]{\@secondoftwo}%
\providecommand \href [0]{\begingroup \@sanitize@url \@href}%
\providecommand \@href[1]{\@@startlink{#1}\@@href}%
\providecommand \@@href[1]{\endgroup#1\@@endlink}%
\providecommand \@sanitize@url [0]{\catcode `\\12\catcode `\$12\catcode
  `\&12\catcode `\#12\catcode `\^12\catcode `\_12\catcode `\%12\relax}%
\providecommand \@@startlink[1]{}%
\providecommand \@@endlink[0]{}%
\providecommand \url  [0]{\begingroup\@sanitize@url \@url }%
\providecommand \@url [1]{\endgroup\@href {#1}{\urlprefix }}%
\providecommand \urlprefix  [0]{URL }%
\providecommand \Eprint [0]{\href }%
\providecommand \doibase [0]{http://dx.doi.org/}%
\providecommand \selectlanguage [0]{\@gobble}%
\providecommand \bibinfo  [0]{\@secondoftwo}%
\providecommand \bibfield  [0]{\@secondoftwo}%
\providecommand \translation [1]{[#1]}%
\providecommand \BibitemOpen [0]{}%
\providecommand \bibitemStop [0]{}%
\providecommand \bibitemNoStop [0]{.\EOS\space}%
\providecommand \EOS [0]{\spacefactor3000\relax}%
\providecommand \BibitemShut  [1]{\csname bibitem#1\endcsname}%
\let\auto@bib@innerbib\@empty
\bibitem [{Note1()}]{Note1}%
  \BibitemOpen
  \bibinfo {note} {We recall that a $G^1-$curve is a tangent vector
  geometrically continuous curve characterized by the fact that the unit
  tangent vector to the curve is continuous \cite {Knott}.}\BibitemShut {Stop}%
\bibitem [{Note2()}]{Note2}%
  \BibitemOpen
  \bibinfo {note} {Nonstandard methods like that based on Particle Swarm
  Optimization proposed in \cite {QuHe} could probably also be applied
  successfully.}\BibitemShut {Stop}%
\bibitem [{Note3()}]{Note3}%
  \BibitemOpen
  \bibinfo {note} {We do not discuss here the particular case in which
  $d_i^{*\prime }$ is also equal to zero, which can be easily
  treated.}\BibitemShut {Stop}%
\end{thebibliography}%


\begin{thebibliography}{99}

\bibitem{wasserman} Wasserman, S. A., Cozzarelli, N. R.: {
Biochemical topology: applications to DNA recombination and replication}.
 Science 232, 951-960 (1986)

\bibitem{delbruck} Delbruck, M.: {Mathematical Problems in the
  Biological Sciences}. Proc. Symp. Appl. Marh. 14, Providence, 
RI: American Mathematical Society, 55 (1962)

\bibitem{vortexLC}  Bowick, M. J.,  Chander, L.,  Schiff, E. A.,
 Srivastava, A. M.: The Cosmological Kibble Mechanism in the Laboratory: String Formation in Liquid Crystals. Science 263, 943-945 (1994)

\bibitem{vortexHe3}  B\"auerle, C., Bunkov, Yu. M., Fisher, S. N.,  Godfrin, H., Pickett, G. R.: Laboratory simulation of cosmic string formation in the early Universe using superfluid $^3$He. Nature 382, 332-334 (1996);
Ruutu, V. M. H., Eltsov, V. B., Gill, A. J., Kibble, T. W. B.,  Krusius,   M., Makhlin, Yu. G., Plaçais, B., Volovik, G. E., Xu, W.: Vortex formation in neutron-irradiated superfluid $^3$He as an analogue of cosmological defect formation. Nature 382, 334-336 (1996)

\bibitem{defectsNColloids} Araki, T., Tanaka, H.: Colloidal Aggregation in a Nematic Liquid Crystal: Topological Arrest of Particles by a Single-Stroke Disclination Line. Phys. Rev. Lett. 97, 127801-4 (2006); 
Tkalec, U., Ravnik, M., \v{C}opar, S., \v{Z}umer, S., Mu\v{s}evi\v{c}, I. : Reconfigurable Knots and Links in Chiral Nematic Colloids. 
Science 333, 62-65 (2011);
\v{Z}umer, S.: 21st International Liquid
Crystal Conference, Keystone, Colorado, USA, July 2–7, 2006,
Kinsley \& Associates, Littletown, (2006) 

\bibitem{Alexander} Alexander, J. W.: Topological Invariants of Knots
and Links. Trans. Amer. Math. Soc. 30, 275-306 (1928)

\bibitem{HOMFLY}
Freyd, P.,  Yetter, D.,  Hoste, J., Lickorish, W. B. R., Millet, K.,
Ocneanu, A.: A new polynomial invariant of knots and links. Bull. AMS 12, 239-246 (1985)

\bibitem{mci}  Kontsevich, M.: Vassiliev’s knot invariants,
Preprint, Max-Planck-Institut f\"ur Mathematik, Bonn;
Bar-Natan, D.: On the Vassiliev knot invariants, Harvard preprint
(1992);
Bar-Natan, D.: On the Vassiliev Knot Invariants, Topology 34, 423-472
(1995)

\bibitem{LR}
Labastida, J. M. F., Ramallo, A. V.: Operator formalism for
Chern-Simons theories.  Phys. Lett. B  227, 92-102 (1989) 

\bibitem{GMM} Guadagnini, E., Martellini, M., Mintchev, M.: Braids and
  quantum group symmetry in Chern-Simons theory. Nucl. Phys. B  336,
  581-609 (1990) 

\bibitem{vologodskii}  Vologodskii, A. V.,  Lukashin, A. V.,
Frank-Kamenetskii, M. D.,  Anshelevich, V. V.: The knot problem in
statistical mechanics of polymer chains.  Zh.  
Eksp. Teor. Fiz.  66, 2153-2163 (1974)

\bibitem{reginald} Nonweiler T. T. F.:
The Numerical Evaluation of Curvilinear Integrals and Areas Defined by
Discrete Data. Inter-university/research councils series (1972)

\bibitem{davis} Davis, P. and Rabinowitz, P.: Methods of Numerical
  Integration.
2nd Ed, New York, Academic Press (1984)

\bibitem{atkinson} Atkinson, K. and Venturino, E.: Numerical
  Evaluation of Line Integrals. Siam J. Numer. Anal. 30, 882-888 (1993)

\bibitem{Knott} Knott, G. D.: Interpolating Cubic Splines. Boston,
  Birkh\"auser (2000)

\bibitem{Dunin} Dunin-Barkowski, P., Sleptsov, A., Smirnov, A.:
Kontsevich integral for knots and Vassiliev invariants.
Int. J. Mod. Phys. A 28, 1330025-62 (2013)

\bibitem{YZFFJSTAT2013} Zhao, Y., Ferrari, F.: A study of polymer
  knots using a simple knot invariant consisting of multiple contour
  integrals. 
 JSTAT 2013, P10010 (2013)


\bibitem{traditional} Davis, P. J., Rabinowitz, P.: Methods of
  numerical integration. Boston, MA, Academic Press, 1984

\bibitem{pivot} Madras, N., Orlitsky, A., Sepp, L. A.: Monte Carlo Generation of Self-Avoiding Walks with Fixed Endpoints and Fixed Length. J. Stat. Phys. 38, 159-183 (1990)


\bibitem{pull} Lesh, N., Mitzenmacher, M., Whitesides, S.: A Complete and
  Effective Move Set for 
Simplified Protein Folding. Proceedings of the Seventh
Annual International Conference on 
Research in Computational Molecular Biology (RECOMB’03), 188-195 (2003)

\bibitem{BFACF} Aragao de Carvalho, C., Caracciolo, S.,
Fr\"ohlich, J.: Polymers and $g|\varphi|^4$ theory in four dimensions.
Nucl. Phys. B 215, 209-248 (1983);
Berg, B., Foerster, D.: Random paths and random surfaces on a digital computer. Phys. Lett. B 106, 323-326 (1981) 

\bibitem{Witten} Witten, E.: Quantum field theory and the Jones
  polynomial. Commun. Math. Phys. 121, 351-399 (1989) 


\bibitem{cotta} Cattaneo, A. S., Cotta-Ramusino, P., Martellini, M.:
  Three-dimensional BF theories 
and the Alexander-Conway invariant of knots.
Nucl. Phys. B 436(1-2), 355-382 (1995)


\bibitem{polyakviro} Polyak, M., Viro, O.: On the Casson knot invariant,
J. Knot Theory Ramifications 10, 711-738 (2001)

\bibitem{QuHe} Qu, L. and He, D.: Solving Numerical Integration
  by Particle Swarm Optimization. 
Work published in the Proceedings of ICICA 2010, Part II,
Zhu, R. et Al. (Eds), Berlin, Heidelberg, Springer Verlag, 228-235
(2010) 


 \bibitem{KoniarisMuthukumar}
 Koniaris, K., Muthukumar, M.: Knottedness in ring polymers. Phys. Rev. 
Lett. 66, 2211-2214 (1991) 


\bibitem{TaylorAsrodi}
Taylor, W. R., Asz\'odi, A.: Protein Geometry, Classification, Topology and Symmetry: A Computational Analysis of Structure. New York, Taylor \& Francis Group, (2005)

\bibitem{bsplines} He, Y., Gu X. F., Qin H.: Automatic Shape Control
  of Triangular B-Splines of Arbitrary Topology. Jour. Comput. Sci.
  \& Technol. 21, 232-237 (2006)


%

%
%
%
%
%
%



\end{thebibliography}
\end{document}